\input epsf

\magnification=1000

\font\titre=cmbx10 scaled\magstep1
\font\auteur=cmr10 scaled\magstep1

\def\square{{\vcenter{\hrule height.4pt
      \hbox{\vrule width.4pt height5pt \hskip5pt
	   \vrule width.4pt}
      \hrule height.4pt}}}
\def\qed{\hfill$\square$}

\def\BN{{\bf N}}
\def\pprod{{\rm prod}}
\def\a{{\alpha}}
\def\b{{\beta}}
\def\c{{\gamma}}
\def\BR{{\bf R}}
\def\Gl{{\rm Gl}}
\def\dpt{{\rm dp}}
\def\BZ{{\bf Z}}
\def\EE{{\cal E}}
\def\BK{{\bf K}}
\def\BQ{{\bf Q}}
\def\f{{\varphi}}
\def\End{{\rm End}}
\def\Id{{\rm Id}} 
\def\HH{{\cal H}}
\def\LL{{\cal L}}
\def\Supp{{\rm Supp}}
\def\CC{{\cal C}}
\def\UU{{\cal U}}


\centerline{\titre Artin monoids inject in their groups}

\bigskip
\centerline{\auteur Luis Paris}
\bigskip
\centerline{\auteur January 31, 2001}

\bigskip\bigskip
\centerline{\bf Abstract}

\bigskip
{\leftskip 1.5cm\rightskip 1.5 cm 
We prove that the natural homomorphism from an Artin monoid to its associated Artin group is always injective.
\par}


\bigskip\bigskip\noindent
{\titre 1. Introduction}

\bigskip\noindent
Let $S$ be a finite set. A {\it Coxeter matrix} over $S$ is a 
matrix $M=(m_{s,t})_{s,t\in S}$ indexed by the elements of 
$S$ and such that:

\smallskip
$\bullet$ $m_{s,s}=1$ for all $s\in S$;

\smallskip
$\bullet$ $m_{s,t}= m_{t,s} \in \{2, 3, 4, \dots, +\infty\}$ 
for all $s,t\in S, s\neq t$.

\smallskip\noindent
A Coxeter matrix $M=(m_{s,t})_{s,t\in S}$ is usually 
represented by its {\it Coxeter graph} $\Gamma$. This is 
defined by the following data:

\smallskip
$\bullet$ S is the set of vertices of $\Gamma$;

\smallskip
$\bullet$ two vertices $s,t \in S$ are joined by an edge if 
$m_{s,t} \ge 3$;

\smallskip
$\bullet$ the edge which joines $s$ and $t$ is labelled by 
$m_{s,t}$ if $m_{s,t}\ge 4$.

\smallskip
The {\it Coxeter system} associated with $\Gamma$ is the 
pair $(W, S)$, where $W$ is the group presented 
by
$$
W = \langle S\ |\ s^2=1\ {\rm for}\ s\in S,\ 
(st)^{m_{s,t}}=1\ {\rm for}\ s,t\in S, s\neq t,\  
m_{s,t}< +\infty \rangle.
$$
The group $W$ is called the {\it Coxeter group} 
associated with $\Gamma$.

Let $\Sigma= \{ \sigma_s; s\in S\}$ be an abstract set in 
one-to-one correspondence with $S$. For two objects $a,b$ 
and $m\in \BN$ we write
$$
\pprod(a,b;m) = \left\{ \matrix{ 
(ab)^{{m \over2}}\hfill&\quad {\rm if}\ m\ {\rm is\ even}\hfill\cr
(ab)^{{m-1 \over 2}}a \hfill&\quad {\rm if}\ m\ {\rm is\ odd}\hfill\cr}
\right.
$$
The {\it Artin system} associated with $\Gamma$ is the pair 
$(G_\Gamma, \Sigma)$, where $G_\Gamma$ is the group 
presented by
$$
G_\Gamma= \langle \Sigma\ |\ \pprod(\sigma_s, \sigma_t ; 
m_{s,t}) = \pprod (\sigma_t, \sigma_s; m_{s,t})\ {\rm for}\ 
s,t\in S, s\neq t, \  m_{s,t}< +\infty \rangle.
$$
The group $G_\Gamma$ is called the {\it Artin group} 
associated with $\Gamma$.

Recall that a {\it monoid} is a semigroup with a unity, and 
a {\it homomorphism} of monoids is a map $\phi: M \to M'$ 
which satisfies $\phi(fg) = \phi(f) \phi(g)$ for all $f,g 
\in M$, and $\phi(1) =1$. The {\it Artin monoid} associated 
with $\Gamma$ is the monoid $G_\Gamma^+$ 
presented by the same generators and relations as $G_\Gamma$.
Let $\iota: G_\Gamma^+ \to 
G_\Gamma$ denote the canonical homomorphism from 
$G_\Gamma^+$ to $G_\Gamma$. The goal of this paper is to 
prove the following.

\bigskip\noindent
{\bf Theorem 1.1.} {\it The homomorphism $\iota: G_\Gamma^+ 
\to G_\Gamma$ is injective for all Coxeter graphs.}

\bigskip
It seems that the authorship of the Artin groups, also 
called generalized braid groups, has to be attributed to 
Jacques Tits, in spite of the fact that his name does not 
always appear in the references. Furthermore, it is in a 
paper of him [Tit2] where these groups appeared for the 
first time. However, it was Brieskorn and Saito who proposed 
in [BS] the question of the study of all these groups 
(Deligne's paper [Del], which appeared at the same time, is 
concerned only with spherical type Artin groups, namely, 
those Artin groups for which $W$ is finite). Some families 
of Artin groups are well understood, but, since the paper of 
Brieskorn and Saito in 1972, very few results concerning all 
Artin groups have been published. In particular, Theorem 1.1 
above was known only for some particular classes, namely, 
for the spherical type Artin groups (see [BS] and [Del]), 
for the two-dimensional Artin groups (see [ChP] and [Cha]), 
and for the FC-type Artin groups (see [Alt] and [Cha]). It 
was unknown, for example, for the so-called affine type Artin groups.

Our proof of Theorem 1.1 is independent of the previous 
approachs of the problem. Note first that, in order to prove 
Theorem 1.1, it suffices to show that there exists an 
injective homomorphism $\psi: G_\Gamma^+ \to G$, where $G$ 
is a group, not necessarily equal to $G_\Gamma$.

We say that a Coxeter graph $\Gamma$ is of {\it small type} 
if $m_{s,t} \in \{2,3\}$ for all $s,t \in S$, $s\neq t$. We 
say that $\Gamma$ has {\it no triangle} if there is no 
triple $\{s,t,r\}$ in $S$ such that $m_{s,t}, m_{s,r}, 
m_{t,r} \ge 3$. The first ingredient in our proof is to show 
that, for any Coxeter graph $\Gamma$, there exists an 
injective homomorphism $\phi: G_\Gamma^+ \to G_{\tilde 
\Gamma}^+$, where $G_{\tilde \Gamma}^+$ is an Artin monoid 
associated to a Coxeter graph $\tilde \Gamma$ of small type 
with no triangle. The homomorphism $\phi$ is obtained by a 
``folding'' as described in [Cri], its construction is 
essentially the same as the one given in [CrP, Sec. 6], and 
the proof of the injectivity is a direct application of 
[Cri, Thm. 1.3]. This construction is given in Section 5.

So, in order to prove Theorem 1.1, it suffices to consider 
only Coxeter graphs of small type with no triangle. Take such 
a Coxeter graph, $\Gamma$. We construct in Section 3 a 
homomorphism $\psi: G_\Gamma^+ \to \Gl (V)$, where $V$ is a 
(infinite dimensional) vector space over $\BQ (x,y)$, and we 
prove in Section 4 that $\psi$ is injective.

If $\Gamma=A_n$, then $G_{A_n}$ is the braid group on $n+1$ 
strings, and $\psi: G_{A_n}^+ \to \Gl(V)$ is equivalent to 
the representation constructed by  Bigelow and Krammer in 
[Big], [Kra1] and 
[Kra2]. In this case, $V$ has finite dimension, and the 
injectivity of $\psi$ implies the injectivity of the induced 
representation $G_{A_n} \to \Gl(V)$. More generally, if 
$\Gamma$ is of spherical and small type, then the 
representation $\psi: G_\Gamma^+ \to \Gl(V)$ is equivalent 
to the ones constructed independently by Digne [Dig], and by 
Cohen and Wales [CW]. In this case again, $V$  has finite 
dimension and the induced representation $G_\Gamma \to 
\Gl(V)$ is injective. We do not know whether the 
representation $G_\Gamma \to \Gl (V)$ induced by $\psi$ is 
injective for all Coxeter graphs of small type with no 
triangle. The construction of $\psi$ and the proof of the 
injectivity are based on a (non always easy) generalization 
of the methods of Krammer, Digne, Cohen and Wales.

\bigskip\noindent
{\bf Acknoledgements.} I would like to thank John Crisp for 
many useful conversations during the preparation of this 
work, and for drawing my attention to the results of [Cri] 
which are one of the main tools of the proof of Theorem 1.1.

\bigskip\bigskip\noindent
{\titre 2. Preliminaries}

\bigskip\noindent
We summarize in this section some well known results on 
Artin monoids, Coxeter groups and root systems, 
and give definitions and some basic properties
of closed sets.
The closed sets have been introduced by Krammer in 
[Kra2] for Artin groups of type $A_n$. 
This notion has been extended 
to the Artin groups of small and spherical type by Digne 
[Dig], Cohen and Wales [CW]. 
Here we extend it to all small type Artin groups.

Let $\Gamma$ be a Coxeter graph. It is shown in [BS] that 
the Artin monoid $G_\Gamma^+$ is {\it cancellative}, namely, 
if $fg_1h = fg_2h$, then $g_1 = g_2$. We say that $h$ is a 
{\it multiple} of $g$ and write $g<h$ if there exists $f\in 
G_\Gamma^+$ such that $gf=h$. The relation $<$ is a partial 
ordering on $G_\Gamma^+$.

Let $\theta: G_\Gamma^+ \to W$ be the homomorphism 
which sends $\sigma_s$ to $s$ for all $s\in S$. Then 
$\theta$ has a natural set-section $\tau: W \to 
G_\Gamma^+$ defined as follows. Let $w\in W$. We 
choose a reduced expression $w=s_1 \dots s_l$ for $w$ and we 
set $\tau(w) = \sigma_{s_1} \dots \sigma_{s_l}$. By Tits' 
solution of the word problem for Coxeter groups [Tit1], the 
definition of $\tau(w)$ does not depend on the choice of the 
reduced expression.

Let $l: W \to \BN$ and $l: G_\Gamma^+ \to \BN$ denote 
the word length functions of $W$ and $G_\Gamma^+$ 
with respect to $S$ and $\Sigma$, respectively. Define a 
partial ordering on $W$ by $u<v$ if $l(v)= l(u)+ l(u^{-1}v)$. 
Then $l(\tau(w)) = l(w)$ for all $w\in W$, and one has $u<v$ 
if and only if $\tau(u) < \tau(v)$.

The proof of the following proposition is essentially the 
same as the one of [Del, Pro. 1.14] and [Mic, Lem. 1.4].

\bigskip\noindent
{\bf Proposition 2.1.} {\it Let $E$ be a finite subset of 
$W$ such that:

\smallskip
$\bullet$ if $u<v$ and $v\in E$, then $u\in E$;

\smallskip
$\bullet$ if $v\in W$ and $s,t \in S$ are such that 
$l(vs) = l(vt) = l(v)+1$ and $vs, vt \in E$, then $m_{s,t} 
<+\infty$ and $v \cdot \pprod( s,t; m_{s,t}) \in E$.

\smallskip\noindent
Then there exists $w_0 \in W$ such that $E=\{ v\in 
W; v<w_0\}$.}

\bigskip
The next proposition is part of [Mic, Prop. 2.1]. It is 
also a direct consequence of [BS, Lem. 2.1 and Prop. 2.3].

\bigskip\noindent
{\bf Proposition 2.2.} {\it Let $f\in G_\Gamma^+$ and let 
$E=\{ w\in W; \tau(w) < f\}$. Then $E$ satisfies:

\smallskip
$\bullet$ if $u<v$ and $v\in E$, then $u\in E$;

\smallskip
$\bullet$ if $v\in W$ and $s,t\in S$ are such that 
$l(vs) = l(vt)= l(v)+1$ and $vs,vt\in E$, then 
$m_{s,t} < +\infty$ and $v \cdot 
\pprod (s,t; m_{s,t}) \in E$.}

\bigskip\noindent
{\bf Definition.} Let $f\in G_\Gamma^+$. By Propositions 2.1 
and 2.2, there exists a unique $w_0 \in W$ such that $\{v\in 
W; \tau(v) <f\} = \{ v\in W; v<w_0\}$. We set
$$
L(f)=w_0.
$$

The next proposition is also part of [Mic, Prop. 2.1]

\bigskip\noindent
{\bf Proposition 2.3.} {\it Let $f,g\in G_\Gamma^+$. Then
$$
L(fg) = L(f \cdot (\tau \circ L)(g)).
$$}

Let $\Pi=\{ \a_s; s\in S\}$ be an abstract set in one-to-one 
correspondence with $S$. The elements of $\Pi$ are called 
{\it simple roots}. Let $U$ denote the real vector space 
having $\Pi$ as a basis, and let $\langle~,~\rangle: U 
\times U \to \BR$ be the symmetric bilinear form on $U$ defined by
$$
\langle \a_s, \a_t \rangle = \left\{ \matrix{
-2 \cos( \pi/ m_{s,t})\hfill& \quad {\rm if}\ m_{s,t} < 
+\infty\hfill\cr
-2\hfill& \quad {\rm if}\ m_{s,t}= +\infty\hfill \cr}
\right.
$$
There is a faithful representation $W \to \Gl(U)$ which is 
defined by
$$
s(x)= x- \langle \a_s, x \rangle \a_s, \quad x\in U,\ s\in S,
$$
and which preserves the bilinear form $\langle~,~\rangle$. 
This representation is called the {\it canonical 
representation} of $W$.

The set $\Phi=\{w\a_s; s\in S, w\in W\}$ is called the {\it 
root system} of $W$. The subsets $\Phi^+ = \{ \sum_{s\in S} 
\lambda_s \a_s \in \Phi; \lambda_s \ge 0$ for all $s\in S\}$ 
and $\Phi^- = \{ \b \in \Phi; -\b \in \Phi^+ \}$ are the 
sets of {\it positive roots} and {\it negative roots}, 
respectively. For $w\in W$ we set $\Phi_w= \{ \b \in \Phi^+; 
w^{-1}\b \in \Phi^- \}$.

We list in the following proposition some well known results 
on root systems (see [Hil] and [Deo]).

\bigskip\noindent
{\bf Proposition 2.4.} {\it (1) $\Phi = \Phi^+ \sqcup \Phi^-$.

\smallskip
(2) $|\Phi_w| = l(w)$ for all $w\in W$.

\smallskip
(3) For all $u,v\in W$ such that $u<v$, one has $\Phi_v = 
\Phi_u \sqcup u\cdot \Phi_{u^{-1} v}$.

\smallskip
(4) For all $w\in W$ and $s\in S$,
$$
l(sw)= \left\{ \matrix{
l(w)+1& \quad {\it if}\ w^{-1} \a_s \in \Phi^+\cr
l(w)-1& \quad {\it if}\ w^{-1} \a_s \in \Phi^- \cr}
\right.
$$

(5) Let $\b =w \a_s \in \Phi^+$, and let $r_\b = wsw^{-1}$. 
Then $r_\b$ acts on $U$ by
$$
r_\b (x) = x- \langle x,\b \rangle \b, \quad x\in U.
$$}

Let $\b\in \Phi^+$. Define the {\it depth} of $\b$ to be
$$
\dpt (\b) = \min\{ l\in \BN;{\rm there\ exists}\ w\in 
W\ {\rm such\ that}\ w \b \in \Phi^-\ {\rm and}\ 
l(w)=l \}.
$$

\noindent
{\bf Lemma 2.5.} {\it Let $\b \in \Phi^+$. Then
$$
\dpt (\b)= \min\{ l\in \BN; {\it there\ exist}\ w\in 
W\ {\it and}\ s\in S\ {\it such\ that}\ \b =w^{-1} 
\a_s\ {\it and}\ l=l(w) +1 \}.
$$}

\noindent
{\bf Proof.} Let $d_1= \min \{ l\in \BN$; there exists $w\in 
W$ such that $w
\b\in \Phi^-$ and $l(w)=l\}$ and $d_2=\min \{ l\in \BN;$ 
there exist $w\in 
W$ and $s\in S$ such that $\b = w^{-1} \a_s$ and 
$l=l(w)+1\}$.

Let $w\in W$ and $s\in S$ such that $\b = w^{-1} \a_s$ and 
$l(w)= d_
2-1$. Since $\b \in \Phi^+$, by Proposition 2.4, $l(sw) =l(w)+1 =d_2$. 
Moreover, $sw\b =s\a_s =-\a_s \in \Phi^-$. This shows that $d_2 \le d_1$.

Let $w\in W$ such that $w\b \in \Phi^-$ and $l(w)=d_1$. Let $s\in S
$ such that $l(sw)= l(w)-1$. Let $v=sw$ and $\c = v\b$. By the 
minimality of $l(w)=d_1$, one has $\c \in \Phi^+$. Moreover, $s\c =w\b \in 
\Phi^-$, thus $\c =\a_s$ and $\b =v^{-1}\a_s$. This shows that $d_1 \le d_2$. 
\qed

\bigskip
The following proposition is proved in [BH, Lem. 1.7].

\bigskip\noindent
{\bf Proposition 2.6.} {\it Let $s\in S$ and $\b \in \Phi^+ \setminus \{ 
\a_s\}$. Then
$$
\dpt (s\cdot \b) = \left\{ \matrix{
\dpt(\b)-1 \hfill & \quad {\it if}\ \langle \a_s, \b \rangle >0 \hfill\cr
\dpt(\b) \hfill & \quad {\it if}\ \langle \a_s, \b \rangle =0 \hfill \cr
\dpt(\b) +1 \hfill & \quad{\it if}\ \langle \a_s, \b \rangle <0 \hfill\cr}
\right.
$$} 

From now on and till the end of the section, we assume that $\Gamma$ is a 
Coxeter graph of small type, namely, that $m_{s,t} \in \{2,3\}$ for all $s,t 
\in S$, $s\neq t$. Note that, under this assumption, all the roots can be 
written $\b = \sum_{s \in S} \lambda_s \a_s$, with $\lambda_s \in \BZ$, and 
one has $\langle \b, \c \rangle \in \BZ$ for all $\b , \c \in \Phi$.

\bigskip\noindent
{\bf Definition.} A subset $
A \subset \Phi^+$ is a {\it closed subset} if:

\smallskip
$\bullet$ $A$ is finite;

\smallskip
$\bullet$ if $\a, \b \in A$, then $\langle \a, \b \rangle \ge -1$;

\smallskip
$\bullet$ if $\a, \b \in A$ and $\langle \a, \b \rangle = -1$, then $\a + 
\b = r_\a(\b) = r_\b (\a) \in A$.

\bigskip\noindent
{\bf Lemma 2.7.} {\it Let $w \in W$. Then $\Phi_w$ is a closed 
subset.}

\bigskip\noindent
{\bf Proof.} Let $\a , \b \in \Phi^+$. A direct calculation shows that: if 
$\langle \a, \b \rangle \le -2$, then $(r_\a r_\b)^l (\a)$ is a positive 
root of the form $a_l \a + b_l \b$, where $a_l, b_l \ge 0$, for all $l\in 
\BN$, and $(r_\a r_\b)^l (\a) \neq (r_\a r_\b)^k (\a)$ for $l \neq k$. This 
implies that: if $\langle \a, \b \rangle \le -2$, then there are infinitely 
many positive roots of the form $a \a + b \b$, with $a,b \ge 0$.

The set $\Phi_w$ is finite since $|\Phi_w| = l(w)$. Let $\a, \b \in \Phi_w
$. If $\c = a \a + b \b$, with $a,b \ge 0$, is a positive root, then $\c \in 
\Phi_w$, since $w^{-1} \c = a w^{-1} \a + b w^{-1} \b$ is a negative root. 
By the above considerations, this implies that 
$\langle \a, \b \rangle \ge -1$ and 
that $\a + \b \in \Phi_w$ if $\langle \a, \b \rangle = -1$.
\qed

\bigskip\noindent
{\bf Proposition 2.8.} {\it Let $A$ be a closed subset of $\Phi^+$ and 
let $E= \{ w\in W; \Phi_w \subset A \}$. Then $E$ satisfies:

\smallskip
$\bullet$ $E$ is finite;

\smallskip
$\bullet$ if $u<v$ and $v\in E$, then $u\in E$;

\smallskip
$\bullet$ if $v\in W$ and $s,t\in S$ are such that $l(vs)= l(vt)= l
(v)+1$ and $vs, vt \in E$, then $v \cdot \pprod (s,t; m_{s,t}) \in E$.}

\bigskip\noindent
{\bf Proof.} If $\Phi_w \subset A$, then $l(w)= |\Phi_w| \le |A|$. Since $A
$ is finite, it follows that $l(w)$ is bounded for all $w\in E$, thus $E$ is 
finite.

Suppose $u<v$ and $v\in E$. Then, by Proposition 2.4, $\Phi_u \subset 
\Phi_v \subset A$, thus $u\in E$.

Let $v \in W$ and $s,t \in S$ such that $l(vs)= l(vt)= l(v)+1$ and 
$vs, vt \in E$. By Proposition 2.4, one has $\Phi_{vs}= \Phi_v \cup \{ v 
\a_s \}$ and $\Phi_{vt} = \Phi_v \cup \{ v \a_t \}$. Let $w = v \cdot 
\pprod (s,t; m_{s,t})$. If $m_{s,t}=2$ then $\Phi_w= \Phi_v \cup \{ v \a_s, 
v\a_t \} \subset A$, thus $w \in E$. If $m_{s,t}=3$, then $\langle v\a_s, 
v\a_t \rangle = \langle \a_s, \a_t \rangle=-1$, thus $v\a_s + v\a_t = v
(\a_s +\a_t) \in A$. It follows that $\Phi_w = \Phi_v \cup \{ v\a_s, v\a_t, 
v(\a_s+\a_t) \} \subset A$, thus $w\in E$.
\qed

\bigskip\noindent
{\bf Definition.} Let $A$ be a closed subset of $\Phi^+$. By 
Propositions 2.1 and 2.8, there exists a unique $w_0 \in W$ such that $\{ w
\in W; \Phi_w \subset A \} = \{ w \in W; w<w_0 \}$. We set
$$
C(A)=w_0.
$$
Note that $C( \Phi_w) = w$ for all $w\in W$.


\bigskip\bigskip\noindent
{\titre 3. The representation}

\bigskip\noindent
Throughout this section, $\Gamma$ is assumed to be a Coxeter graph of small 
type with no triangle, namely, $m_{s,t}\in \{2,3\}$ for all $s,t \in S$,
$s \neq t$, and there is no triple $\{s,t,r\}$ in $S$ such that
$m_{s,t} = m_{s,r} = m_{t,r} =3$. 
Our aim here is to construct a (infinite 
dimensional) linear representation $\psi: G_\Gamma^+ \to \Gl (V)$. We will 
prove in Section 4 that this linear representation is faithful. This will 
imply that $\iota : G_\Gamma^+ \to G_\Gamma$ is injective. 

Let $\EE= \{ e_\b; \b \in \Phi^+ \}$ be an abstract set in one-to-one 
correspondence with $\Phi^+$, let $\BK = \BQ (x,y)$ denote the field of 
rationnal functions on two variables over $\BQ$, and let $V$ be the 
$\BK$-vector space having $\EE$ as a basis.

For all $s\in S$, we define a linear transformation $\f_s : V \to V$ by
$$
\varphi_s (e_\b) = \left\{ \matrix{
0 \hfill &\quad{\rm if}\ \b=\a_s\hfill\cr
e_\b \hfill& \quad {\rm if}\ \langle \a_s,\b \rangle =0 \hfill\cr
y \cdot e_{\b-a\a_s} \hfill& \quad {\rm if}\ \langle \a_s,\b \rangle =a>0\ 
{\rm and}\ \b \neq \a_s\hfill\cr
(1-y) \cdot e_\b + e_{\b+a\a_s}\hfill & \quad {\rm if}\ \langle \a_s,\b 
\rangle = -a<0 \hfill\cr}
\right.
$$
A direct (case by case) calculation shows that
\bigskip
\centerline{\vbox{\halign{\hfill$#$\ &$=\ #$\hfill&\quad$#$\hfill\cr
\f_s \f_t&\f_t \f_s&{\rm if}\ m_{s,t}=2\cr
\f_s \f_t \f_s& \f_t \f_s \f_t&{\rm if}\ m_{s,t}=3\cr}}}
\bigskip\noindent
So: 

\bigskip\noindent
{\bf Proposition 3.1.} {\it The mapping $\sigma_s \to \f_s$, $s \in S$, 
induces a homomorphism $\f : G_\Gamma^+ \to \End (V)$.}

\bigskip
Now, for all $s \in S$ and all $\b \in \Phi^+$, take a polynomial $T(s,\b) 
\in \BQ [y]$ and define $\psi_s: V \to V$ by
$$
\psi_s (e_\b) = \f_s(e_\b) +xT(s,\b) \cdot e_{\a_s}.
$$
The goal of this section is to prove the following:

\bigskip\noindent
{\bf Theorem 3.2.} {\it There is a choice of polynomials $T(s,\b)$, $s  \in 
S$ and $\b \in \Phi^+$, so that the mapping $\sigma_s \to \psi_s$, $s \in S
$, induces a homomorphism $\psi: G_\Gamma^+ \to \Gl(V)$.}

\bigskip
Let $s \in S$ and $\b \in \Phi^+$. We define the polynomial $T(s,\b)$ by 
induction on $\dpt (\b)$. Assume first that $\dpt(\b)=1$. There exists $t 
\in S$ such that $\b = \a_t$. Then we set
\bigskip
\vbox{\halign{#\quad\hfill&\hfill$#$&$#$\hfill&\quad$#$\hfill\cr
(D1)& T(s,\a_t)= & \ y^2&{\rm if}\ t=s\cr
\noalign{\smallskip}
(D2)& T(s,\a_t)= & \ 0&{\rm if}\ t \neq s\cr}}
\bigskip\noindent
Now, assume that $\dpt(\b) \ge 2$. We choose $t \in S$ such that $\dpt(t 
\cdot \b) = \dpt (\b) -1$. By Proposition 2.6, one has $ \langle \a_t, \b 
\rangle =b >0$.

\smallskip\noindent
{\it Case 1:} $\langle \a_s, \b \rangle = a >0$. Then we set
\bigskip
\vbox{\halign{#\quad\hfill&\hfill$#$\ =&\ $#$\hfill\cr
(D3)&T(s,\b)&y^{\dpt(\b)} (y-1)\cr}}
\bigskip

\noindent
{\it Case 2:} $\langle \a_s,\b \rangle = 0$. Then we set
\bigskip
\vbox{\halign{#\quad\hfill&\hfill$#$\ =&\ $#$\hfill&\quad$#$\hfill\cr
(D4)&T(s,\b)& y \cdot T(s, \b-b\a_t) &{\rm if}\ \langle \a_s,\a_t \rangle =
0\cr
\noalign{\smallskip}
(D5)&T(s,\b)& (y-1) \cdot T(s, \b-b\a_t) + y \cdot T(t, \b -b\a_s -b\a_t) &
{\rm if}\ \langle \a_s,\a_t \rangle=-1\cr}}
\bigskip

\noindent
{\it Case 3:} $\langle \a_s,\b \rangle = -a <0$. Then we set
\bigskip
\vbox{\halign{#\quad\hfill&\hfill$#$\ =&\ $#$\hfill&\quad$#$\hfill\cr
(D6) & T(s,\b) &y \cdot T(s, \b -b\a_t) &{\rm if}\ \langle \a_s,\a_t 
\rangle =0\cr
\noalign{\smallskip}
(D7) &T(s,\b) &(y-1) \cdot T(s, \b -b\a_t) +y \cdot T(t, \b -(b-a)\a_s -b
\a_t) &{\rm if}\ \langle \a_s,\a_t \rangle =-1\ {\rm and}\ b>a\cr
\noalign{\smallskip}
(D8) & T(s,\b) & T(t, \b -b\a_t) +(y-1) \cdot T(s, \b -b\a_t) &{\rm if}\ 
\langle \a_s,\a_t \rangle =-1\ {\rm and}\ b=a\cr
\noalign{\smallskip}
(D9) &T(s,\b)& y\cdot T(s, \b -b\a_t) +T(t, \b -b\a_t) + y^{\dpt(\b)-1} (1-
y) &{\rm if}\ \langle \a_s, \a_t \rangle =-1\ {\rm and}\ b<a\cr}}
\bigskip

The proofs of the following lemmas 3.3 and 3.4 are long and tedious case by 
case verifications and they are not very instructive for the remainder of 
the paper. So, we put them in a separate section at the end of the paper 
and continue with the proof of Theorem 3.2.

\bigskip\noindent
{\bf Lemma 3.3.} {\it Let $s\in S$ and $\b \in \Phi^+$ such that $\dpt(\b) 
\ge 2$ and $\langle \a_s,\b \rangle =0$. Then the definition of $T(s, \b)$ 
does not depend on the choice of the $t \in S$ such that $\dpt( t\cdot \b) 
= \dpt(\b) -1$.}

\bigskip\noindent
{\bf Lemma 3.4.} {\it Let $s\in S$ and $\b \in \Phi^+$ such that $\dpt(\b) 
\ge 2$ and $\langle \a_s, \b \rangle = -a<0$. Then the definition of $T(s, 
\b)$ does not depend on the choice of the $t \in S$ such that $\dpt(t \cdot 
\b) = \dpt(\b)-1$.}

\bigskip\noindent
{\bf Lemma 3.5.} {\it Let $s,t \in S$ and $\b \in \Phi^+$ such that 
$\langle \a_s,\a_t \rangle =-1$, $\langle \a_s, \b \rangle =0$, and $\langle 
\a_t, \b \rangle =0$. Then
$$
T(s,\b)=T(t,\b).
$$}

\noindent
{\bf Proof.} We argue by induction on $\dpt(\b)$. 
Assume first that $\dpt(\b) =1
$. There exists $r\in S$ such that $\b = \a_r$. One has $r \neq s$ and $r 
\neq t$ since $\langle \a_s, \b \rangle = \langle \a_t, \b \rangle =0$. 
Then, by (D2),
$$
T(s,\b) = T(s,\a_r) = 0 = T(t,\a_r) = T(t,\b).
$$
Now, assume that $\dpt(\b) \ge 2$. We choose $r\in S$ such that $\dpt(r 
\cdot \b) = \dpt(\b) -1$. By Proposition 2.6, one has $\langle \a_r, \b 
\rangle =c >0$.

\smallskip\noindent
{\it Case 1:} $\langle \a_s, \a_r \rangle =0$ and $\langle \a_t, \a_r 
\rangle =0$. Then
\bigskip
\centerline{\vbox{\halign{\hfill$#$\ &=\ $#$\hfill&\quad #\hfill\cr
T(s,\b) &y \cdot T(s, \b -c\a_r) &by (D4)\cr
&y \cdot T(t,\b-c\a_r)&by induction\cr
&T(t,\b)&by (D4)\cr}}}
\bigskip

\noindent
{\it Case 2:} $\langle \a_s, \a_r \rangle =0$ and $\langle \a_t, \a_r \rangle 
=-1$. We cannot have $\dpt(\b) \ge 3$ in this case. Suppose $\dpt(\b) \ge 4
$. Then
\bigskip
\centerline{\vbox{\halign{\hfill$#$\ &=\ $#$\hfill&\quad #\hfill\cr
T(s,\b) &y \cdot T(s, \b -c\a_r) &by (D4)\cr
& y(y-1) \cdot T(s, \b -c \a_t -c\a_r) +y^2 \cdot T(t, \b -c\a_s -c\a_t -c
\a_r) &by (D5)\cr
&y^{\dpt(\b)-1} (y-1)^2 +y^2 \cdot T(t, \b -c\a_s -c\a_t -c\a_r) &by (D3)
\cr
&y^{\dpt(\b)-1} (y-1)^2 +y^2 \cdot T(r, \b -c\a_s -c\a_t -c\a_r) &by 
induction\cr
& (y-1) \cdot T(t, \b -c\a_r) + y\cdot T(r, \b -c\a_t -c\a_r)& by (D3) and 
(D4)\cr
& T(t,\b) &by (D5).\cr}}} 
\bigskip

Since $\Gamma$ has no triangle, we cannot have $\langle \a_s, \a_r \rangle 
=-1$ and $\langle \a_t, \a_r \rangle =-1$ ( because $\langle \a_s, \a_t 
\rangle =-1$). So, Case 1 and Case 2 are the only possible cases. \qed

\bigskip\noindent
{\bf Lemma 3.6.} {\it Let $s,t \in S$ such that $m_{s,t}=2$. Then $\psi_s 
\psi_t = \psi_t \psi_s$.}

\bigskip\noindent
{\bf Proof.} Let $\b \in \Phi^+$. We compute $(\psi_s \psi_t) (e_\b)$ and 
$(\psi_t \psi_s) (e_\b)$ replacing $T(s,\a_s)$ and $T(t,\a_t)$ by $y^2$, 
and replacing $T(s,\a_t)$ and $T(t, \a_s)$ by $0$, and we compare both 
expressions. This can be easily made with a computer.

\smallskip\noindent
{\it Case 1:} $\b = \a_s$. Then we directly obtain $(\psi_s \psi_t) (e_\b) 
= (\psi_t \psi_s) (e_\b)$.

\smallskip\noindent
{\it Case 2:} $\langle \a_s, \b \rangle =0$ and $\langle \a_t, \b 
\rangle =0$. Then we directly obtain $(\psi_s \psi_t) (e_\b) = (\psi_t 
\psi_s) (e_\b)$.

\smallskip\noindent
{\it Case 3:} $\langle \a_s, \b \rangle = 0$ and $\langle \a_s, \b \rangle 
= b>0$. Then the equality $(\psi_s \psi_t)(e_\b) = (\psi_t \psi_s) (e_\b)$ 
is equivalent to
$$
T(s,\b) = y \cdot T(s, \b -b\a_t).
$$
This equality follows from (D4).

\smallskip\noindent
{\it Case 4:} $\langle \a_s, \b \rangle =0$ and $\langle \a_t, \b \rangle = 
-b <0$. Then the equality $(\psi_s \psi_t)(e_\b) = (\psi_t \psi_s) (e_\b)$ 
is equivalent to
$$
T(s, \b +b\a_t) = y \cdot T(s, \b).
$$
This equality follows from (D4).

\smallskip\noindent
{\it Case 5:} $\langle \a_s, \b \rangle = a>0$ and $\langle \a_t, \b 
\rangle = b>0$. Then the equality $(\psi_s \psi_t)(e_\b) = (\psi_t \psi_s) 
(e_\b)$ is equivalent to
$$\displaylines{
T(s,\b)= y \cdot T(s,\b -b\a_t),\cr
T(t,\b) = y \cdot T(t,\b -a\a_s).\cr}
$$
These two equalities follow from (D3).

\smallskip\noindent
{\it Case 6:} $\langle \a_s, \b \rangle = a >0$ and $\langle \a_t, \b 
\rangle =-b<0$. Then the equality $(\psi_s \psi_t)(e_\b) = (\psi_t \psi_s) 
(e_\b)$ is equivalent to
$$\displaylines{
T(s, \b +b\a_t) = y \cdot T(s,\b),\cr
T(t,\b) = y \cdot T(t, \b -a \a_s).\cr}
$$
The first equality follows from (D3) and the second one from (D6).

\smallskip\noindent
{\it Case 7:} $\langle \a_s, \b \rangle = -a<0$ and $\langle \a_t, \b 
\rangle =-b <0$. Then the equality $(\psi_s \psi_t)(e_\b) = (\psi_t \psi_s) 
(e_\b)$ is equivalent to
$$\displaylines{
T(s,\b +b\a_t)= y \cdot T(s,\b), \cr
T(t, \b +a\a_s) = y \cdot T(t,\b).\cr}
$$
These two equalities follow from (D6). \qed

\bigskip\noindent
{\bf Lemma 3.7.} {\it Let $s,t \in S$ such that $m_{s,t}=3$. Then $\psi_s 
\psi_t \psi_s = \psi_t \psi_s \psi_t$.}

\bigskip\noindent
{\bf Proof.} Let $\b \in \Phi^+$. We conpute $(\psi_s \psi_t \psi_s) (e_\b)
$ and $(\psi_t \psi_s \psi_t) (e_\b)$ replacing $T(s,\a_s)$ and $T(t,\a_t)
$ by $y^2$, replacing $T(s,\a_t)$ and $T(t,\a_s)$ by $0$, and replacing $T
(s, \a_s+\a_t)$ and $T(t, \a_s+\a_t)$ by $y^2(y-1)$, and we compare both 
expressions.

\smallskip\noindent
{\it Case 1:} $\b = \a_s$. Then we directly obtain $(\psi_s \psi_t \psi_s) 
(e_\b) = (\psi_t \psi_s \psi_t) (e_\b)$.

\smallskip\noindent
{\it Case 2:} $\b = \a_s +\a_t$. Then we directly obtain $(\psi_s \psi_t 
\psi_s) (e_\b) = (\psi_t \psi_s \psi_t) (e_\b)$.

\smallskip\noindent
{\it Case 3:} $\langle \a_s, \b \rangle =0$ and $\langle \a_t, \b \rangle =
0$. Then the equality $(\psi_s \psi_t \psi_s) (e_\b) = (\psi_t \psi_s
\psi_t) (e_\b)$ is equivalent to
$$
T(s,\b) = T(t,\b).
$$
This equality follows from Lemma 3.5.

\smallskip\noindent
{\it Case 4:} $\langle \a_s, \b \rangle =0$ and $\langle \a_t, \b \rangle 
=b>0$. Then the equality $(\psi_s \psi_t \psi_s) (e_\b) = (\psi_t \psi_s 
\psi_t) (e_\b)$ is equivalent to
$$\displaylines{
T(t,\b)= y \cdot T(s, \b -b\a_t),\cr
(1-y) \cdot T(t,\b) +y \cdot T(s,\b) = y^2 \cdot T(t, \b -b\a_s -b\a_t).
\cr}
$$
The first equality follows from (D3) and the second one follows from
the first one and from (D5).

\smallskip\noindent
{\it Case 5:} $\langle \a_s, \b \rangle =0$ and $\langle \a_t, \b \rangle 
=-b<0$. Then the equality $(\psi_s \psi_t \psi_s) (e_\b) = (\psi_t \psi_s 
\psi_t) (e_\b)$ is equivalent to
$$\displaylines{
(1-y) \cdot T(s,\b) + T(s, \b +b\a_t) = T(t,\b),\cr
y\cdot T(s,\b) = (1-y) \cdot T(t, \b+b\a_t) + T(t, \b+ b\a_s +b\a_t).\cr}
$$
The first equality follows from (D8) and the second one from (D5).

\smallskip\noindent
{\it Case 6:} $\langle \a_s, \b \rangle = a>0$ and $\langle \a_t, \b 
\rangle =b>0$.  Then the equality $(\psi_s \psi_t \psi_s) (e_\b) = (\psi_t 
\psi_s \psi_t) (e_\b)$ is equivalent to
$$\displaylines{
y \cdot T(s,\b -a\a_s -(a+b)\a_t) = (1-y) \cdot T(s, \b -b\a_t) +T(t, \b),
\cr
(1-y) \cdot T(t, \b -a\a_s) +T(s,\b) = y \cdot T(t, \b -(a+b)\a_s -b\a_t),
\cr
T(t,\b -a\a_s) = T(s, \b-b\a_t).\cr}
$$
These three equalities follow from (D3).

\smallskip\noindent
{\it Case 7:} $\langle \a_s, \b \rangle = a>0$, $\langle \a_t, \b \rangle 
=-b<0$, and $a>b$. Then the equality $(\psi_s \psi_t \psi_s) (e_\b) = 
(\psi_t \psi_s \psi_t) (e_\b)$ is equivalent to
$$\displaylines{
y^2 \cdot T(s, \b -a\a_s -(a-b)\a_t) = (1-y)^2 \cdot T(s,\b) + (1-y) \cdot 
T(s, \b +b\a_t) +y \cdot T(t,\b),\cr
T(s,\b)= T(t, \b -(a-b)\a_s +b\a_t),\cr
y\cdot T(t, \b -a\a_s) = (1-y) \cdot T(s,\b) + T(s, \b+b\a_t).\cr}
$$
The second and third equalities follow from (D3), and the first one follows 
from the third one and from (D7).

\smallskip\noindent
{\it Case 8:} $\langle \a_s, \b \rangle =a >0$, $\langle \a_s, \b \rangle 
=-b<0$, and $a=b$.  Then the equality $(\psi_s \psi_t \psi_s) (e_\b) = 
(\psi_t \psi_s \psi_t) (e_\b)$ is equivalent to
$$\displaylines{
y \cdot T(s, \b -a\a_s) = (1-y)^2 \cdot T(s,\b) + (1-y) \cdot T(s, \b +a
\a_t) +y \cdot T(t,\b), \cr
T(t, \b +a\a_t) = y \cdot T(s, \b),\cr
y \cdot T(t, \b -a\a_s) = (1-y) \cdot T(s,\b) + T(s,\b +a\a_t).\cr}
$$
The second equality follows from (D3), the third one follows from (D5), and 
the first one follows from the third one and from (D8).

\smallskip\noindent
{\it Case 9:} $\langle \a_s, \b \rangle =a>0$, $\langle \a_t, \b \rangle =-b
<0$, and $a<b$.  Then the equality $(\psi_s \psi_t \psi_s) (e_\b) = (\psi_t 
\psi_s \psi_t) (e_\b)$ is equivalent to
$$\displaylines{
y(1-y) \cdot T(s, \b -a\a_s) +y \cdot T(s, \b -a\a_s +(b-a)\a_t) = (1-y)^2 
\cdot T(s,\b) +(1-y) \cdot T(s, \b +b\a_t) +y \cdot T(t,\b),\cr
y \cdot T(s,\b) = (1-y) \cdot T(t, \b +b\a_t) +T(t, \b +(b-a)\a_s +b\a_t),
\cr
y \cdot T(t, \b -a\a_s) = (1-y) \cdot T(s,\b) + T(s, \b +b\a_t).\cr}
$$
The second equality follows from (D3), the third one follows from (D7), and 
the first one follows from the third one and from (D9).
 
\smallskip\noindent
{\it Case 10:} $\langle \a_s, \b \rangle =-a<0$ and $\langle \a_t, \b 
\rangle =-b<0$. Then the equality $(\psi_s \psi_t \psi_s) (e_\b) = (\psi_t 
\psi_s \psi_t) (e_\b)$ is equivalent to
$$\displaylines{
(1-y) \cdot T(s, \b +a\a_s) +T(s, \b +a\a_s +(a+b)\a_t) = y \cdot T(t,\b),
\cr
y \cdot T(s,\b)= (1-y) \cdot T(t, \b+ b\a_t) +T(t, \b +(a+b)\a_s +b\a_t),\cr
(1-y) \cdot T(t,\b) + T(t, \b +a\a_s) = (1-y) \cdot T(s,\b) + T(s, \b +b
\a_t).\cr}
$$
The first and second equalities follow from (D7), and the third one follows 
from (D9). \qed

\bigskip\noindent
{\bf Lemma 3.8.} {\it Let $s \in S$. Then $\psi_s$ is invertible.}

\bigskip\noindent
{\bf Proof.} Let $\rho_s: V \to V$ be the linear transformation defined by
$$
\rho_s(e_\b)= \left\{ \matrix{
x^{-1}y^{-2} \cdot e_{\a_s} \hfill&\quad{\rm if}\ \b = \a_s \hfill\cr
e_\b -y^{-2} T(s,\b) \cdot e_{\a_s}\hfill&\quad {\rm if}\ \langle 
\a_s, \b \rangle =0 \hfill \cr
(1-y^{-1}) \cdot e_\b + e_{\b -a\a_s} -y^{-2} T(s, \b -a\a_s) \cdot e_
{\a_s} \hfill\cr
\quad + y^{-2} (y^{-1}-1) T(s,\b) \cdot e_{\a_s} \hfill&
\quad {\rm if}\ \langle \a_s,\b \rangle =a>0\ {\rm and}\ \b 
\neq \a_s \hfill\cr
y^{-1} \cdot e_{\b +a\a_s} -y^{-3} T(s, \b +a\a_s) \cdot e_{\a_s} 
\hfill & \quad{\rm if}\ \langle \a_s, \b \rangle =-a<0 \hfill\cr}
\right.
$$
A direct case by case calculation shows that
$\psi_s \circ \rho_s = \rho_s \circ \psi_s = 
\Id_V$. So, $\psi_s$ is invertible. \qed

\bigskip
This finishes the proof of Theorem 3.2. 

\bigskip\bigskip\noindent
{\titre 4. Faithfulness}

\bigskip\noindent
Throughout this section, $\Gamma$ is again assumed to be a 
Coxeter graph of small type with no triangle. Our goal here 
is to prove the following.

\bigskip\noindent
{\bf Theorem 4.1.} {\it The representation $\psi: G_\Gamma^+ 
\to \Gl (V)$ defined in Section 3 is faithful.}

\bigskip
Since $\Gl (V)$ is a group, it follows:

\bigskip\noindent
{\bf Corollary 4.2.} {\it The homomorphism $\iota: 
G_\Gamma^+ \to G_\Gamma$ is injective.}

\bigskip
Let $V_+= \oplus_{ \b \in \Phi^+} \BQ [x,y] e_\b$ denote the 
free $\BQ[x,y]$-module having $\EE= \{ e_\b; \b\in \Phi^+\}$ 
as a basis. The coefficients of $\psi(g)$ lie in $\BQ [x,y]$, 
for all $g \in G_\Gamma^+$, thus $V_+$ is invariant by the 
action of $G_\Gamma^+$. We denote by $\psi_+: 
G_\Gamma^+ \to \End(V_+)$ the restriction of $\psi$ to 
$V_+$.

Let $V_0= \oplus_{\b \in \Phi^+} \BR e_\b$ denote the real 
vector space having $\EE$ as a basis. Replacing $x$ by $0$ 
and $y$ by a value $0 <y_0< 1$, the homomorphism $\psi_+$ 
induces a homomorphism $\psi_0: G_\Gamma^+ \to \End( V_0)$.

Let $\HH$ denote the Hilbert space of series $\sum_{\b \in 
\Phi^+} \lambda_\b e_\b$ such that $\sum_{\b \in \Phi^+} 
\lambda_\b^2 <+\infty$, and  let $\LL( \HH)$ denote the 
space of continuous linear transformations of $\HH$. Remark 
that $V_0$ is dense in $\HH$. For all $s\in S$ and all $\b 
\in \Phi^+$, one has
$$
\psi_0(\sigma_s) (e_\b)= \left\{ \matrix{
0\hfill &\quad {\rm if}\ \b = \a_s\hfill\cr
e_\b\hfill &\quad {\rm if}\ \langle \a_s,\b \rangle =0 
\hfill\cr
y_0 \cdot e_{\b -a\a_s} \hfill&\quad {\rm if}\ \langle 
\a_s,\b \rangle =a>0\ {\rm and}\ \b \neq \a_s \hfill\cr
(1-y_0) \cdot e_\b + e_{\b +a\a_s} \hfill &\quad {\rm if}\ 
\langle \a_s,\b \rangle =-a<0\hfill \cr}
\right.
$$
The hypothesis $0<y_0<1$ implies that $\| \psi_0 (\sigma_s) 
(e_\b) \|_2 \le \sqrt{2}$ for all $\b \in \Phi^+$, 
thus $\psi_0 (\sigma_s) \in \LL( 
\HH)$. So, $\psi_0: G_\Gamma^+ \to \End(V_0)$ induces a 
homomorphism $\psi_\infty: G_\Gamma^+ \to \LL( \HH)$.

\bigskip\noindent
{\bf Definition.} Let $A$ be a subset of $\Phi^+$. Then 
$U_A$ denotes the set of series $\sum \lambda_\b e_\b \in 
\HH$ such that:

\smallskip
$\bullet$ $\lambda_\b \ge 0$ for all $\b \in \Phi^+$;

\smallskip
$\bullet$ $\lambda_\b=0$ if and only if $\b \in A$.

\smallskip\noindent
Note that $U_A$ is nonempty, even if $\Phi^+ \setminus A$ is 
infinite, and one has $U_A \cap U_B = \emptyset$ if $A \neq 
B$.

\bigskip\noindent
{\bf Lemma 4.3.} {\it Let $A \subset \Phi^+$ and $g \in 
G_\Gamma^+$. There exists a unique subset $B \subset \Phi^+$ 
such that $\psi_\infty (g) \cdot U_A \subset U_B$.}

\bigskip\noindent
{\bf Proof.} The hypothesis $0<y_0<1$ implies that the 
coefficients of $\psi_\infty (\sigma_s)$ are $\ge 0$ for all 
$s\in S$, thus the coefficients of $\psi_\infty (g)$ are 
$\ge 0$. Let $\psi_\infty (g) (e_\b)= \sum_{\c \in \Phi^+} 
a_{\c,\b} e_\c$, and let $\Supp_g(e_\b)$ denote the set of 
$\c \in \Phi^+$ such that $a_{\c,\b} >0$. Let
$$
A' =\Phi^+ \setminus A, \quad B'= \cup_{\b \in A'} \Supp_g 
(e_\b), \quad B= \Phi^+ \setminus B'.
$$
Then $\psi_\infty (g) \cdot U_A \subset U_B$. \qed

\bigskip\noindent
{\bf Definition.} Let $A \subset \Phi^+$ and $g \in 
G_\Gamma^+$. Then $g \ast A =B$ denotes the unique subset $B 
\subset \Phi^+$ such that $\psi_\infty (g) \cdot U_A \subset 
U_B$.

\bigskip\noindent
{\bf Lemma 4.4.} {\it Let $A \subset \Phi^+$ and $s \in S$. 
Then
\bigskip
\centerline{\vbox{\halign{\hfill$#$&$#$\hfill\cr
\sigma_s \ast A = \{ \a_s \}\ &\cup\ \{ \b \in \Phi^+; 
\langle \a_s, \b \rangle =0\ {\it and}\ \b \in A \} \cr
&\cup\ \{ \b \in \Phi^+; \langle \a_s, \b \rangle =a>0,\ \b 
\neq \a_s,\ {\it and}\ \b-a\a_s \in A \} \cr
&\cup\ \{ \b \in \Phi^+; \langle \a_s, \b \rangle =-a<0,\ \b 
\in A,\ {\it and}\ \b+a\a_s \in A \}.\cr}}}
\bigskip}

\noindent
{\bf Proof.} Let $\b \in \Phi^+$. Then
$$
{\rm Supp}_{\sigma_s} (e_\b) = \left\{ \matrix{
\emptyset \hfill &\quad {\rm if}\ \b=\a_s \hfill \cr
\{e_\b\} \hfill &\quad {\rm if}\ \langle \a_s,\b \rangle =0 
\hfill \cr
\{ e_{\b-a\a_s} \} \hfill &\quad {\rm if}\ \langle \a_s, \b 
\rangle =a>0\ {\rm and}\ \b \neq \a_s\hfill \cr
\{ e_\b, e_{\b+a\a_s} \} \hfill &\quad {\rm if}\ \langle 
\a_s, \b \rangle =-a<0 \hfill \cr}
\right.
$$
Let $A'= \Phi^+ \setminus A$ and $B' =\cup_{\b \in A'} 
\Supp_{\sigma_s} (e_\b)$. Then
\bigskip
\centerline{\vbox{\halign{
\hfill$#$&\hfill$#$&$#$\hfill\cr
B'=&&\{ \b \in \Phi^+; \langle \a_s, \b \rangle =0\ {\rm 
and}\ \b \in A' \}\cr
&\cup\ &\{ \b \in \Phi^+; \langle \a_s, \b \rangle =a>0,\ \b 
\neq \a_s,\ {\rm and}\ \b-a\a_s \in A' \}\cr
&\cup\ &\{ \b \in \Phi^+; \langle \a_s, \b \rangle =-a<0,\ 
{\rm and\ either}\ \b \in A'\ {\rm or}\ \b+a\a_s \in A' 
\}\cr}}}
\bigskip\noindent
thus
\bigskip
\centerline{\vbox{\halign{
\hfill$#$&$#$\hfill&$#$\hfill\cr
\sigma_s \ast A =&\multispan 2 $\ B = \Phi^+ \setminus B'$ 
\hfill \cr
= &\ \{ \a_s\}\ &\cup\ \{ \b \in \Phi^+; \langle \a_s, \b 
\rangle =0\ {\rm and}\ \b \in A\} \cr
&&\cup\ \{ \b \in \Phi^+; \langle \a_s, \b \rangle =a>0,\ \b 
\neq \a_s,\ {\rm and}\ \b-a\a_s \in A\}\cr
&&\cup\ \{\b \in \Phi^+; \langle \a_s, \b \rangle =-a<0,\ \b 
\in A,\ {\rm and}\ \b+a\a_s \in A\}. \quad \square\cr}}}
\bigskip

\noindent
{\bf Remark.} Let $A \subset \Phi^+$ and $s\in S$. Then
\bigskip
\centerline{\vbox{\halign{
\hfill$#$&$#$\hfill&$#$\hfill\cr
s(A \setminus \{ \a_s\}) =&&\{ \b \in \Phi^+; \langle \a_s, 
\b \rangle =0\ {\rm and}\ \b \in A \}\cr
&\cup\ &\{ \b \in \Phi^+; \langle \a_s, \b \rangle =a>0,\ \b 
\neq \a_s,\ {\rm and}\ \b-a\a_s \in A \}\cr
&\cup\ &\{ \b \in \Phi^+; \langle \a_s, \b \rangle =-a<0\ 
{\rm and}\ \b+a\a_s \in A\}.\cr}}}
\bigskip\noindent
In particular, one has
$$
\sigma_s \ast A \subset \{ \a_s \} \cup s( A \setminus \{ 
\a_s\}).
$$

\noindent
{\bf Lemma 4.5.} {\it Let $A$ be a closed subset of $\Phi^+$ 
and $s \in S$. Then $\sigma_s \ast A$ is also a closed 
subset.}

\bigskip\noindent
{\bf Proof.} Since $A$ is finite, $\sigma_s \ast A$ is also 
finite. Now, we take $\b_1, \b_2 \in \sigma_s \ast A$ and we 
prove:

\smallskip
$\bullet$ $\langle \b_1, \b_2 \rangle \ge -1$;

\smallskip
$\bullet$ if $\langle \b_1, \b_2 \rangle =-1$, then 
$\b_1+\b_2 \in \sigma_s \ast A$.

\smallskip
Assume first that $\b_1 = \a_s$. If $\langle \a_s, \b_2 
\rangle =-a<0$, then $\b_2, \b_2+a\a_s \in A$ (by Lemma 
4.4), thus
$$
\langle \b_2, \b_2+a\a_s \rangle = \langle \b_2, \b_2 
\rangle + a \langle \b_2, \a_s \rangle =2-a^2 \ge -1
$$
(since $A$ is closed), therefore $a=1$. This shows that 
$\langle \a_s, \b_2 \rangle \ge -1$. Suppose 
$\langle \a_s, \b_2 \rangle =-1$. Then $\langle \a_s, 
\b_2+\a_s \rangle =1$ and $\b_2 \in A$ (by Lemma 4.4), 
thus, by Lemma 4.4, 
$\b_2+\a_s \in \sigma_s \ast A$.

Assume now that $\b_1 \neq \a_s$ and $\b_2 \neq \a_s$. One 
has $\b_1, \b_2 \in s(A \setminus \{\a_s\})$, thus $s(\b_1), 
s(\b_2) \in A$, therefore $\langle \b_1, \b_2 \rangle = 
\langle s(\b_1), s(\b_2) \rangle \ge -1$ (since $A$ is 
closed). Suppose now that $\langle \b_1, \b_2 \rangle =-1$.

\smallskip\noindent
{\it Case 1:} $\langle \a_s, \b_1 \rangle =0$ and $\langle 
\a_s, \b_2 \rangle =0$. Then $\langle \a_s, \b_1 + \b_2 
\rangle =0$. Moreover, one has $\b_1, \b_2 \in A$ (by Lemma 
4.4), thus $\b_1+\b_2 \in A$ (since $A$ is closed), 
therefore $\b_1 + \b_2 \in \sigma_s \ast A$ (by Lemma 4.4).

\smallskip\noindent
{\it Case 2:} $\langle \a_s, \b_1 \rangle =0$ and $\langle 
\a_s, \b_2 \rangle =b>0$. Then $\langle \a_s, \b_1+\b_2 
\rangle =b>0$. Moreover, one has $\b_1, \b_2-b\a_s \in A$ 
(by Lemma 4.4) and $\langle \b_1, \b_2-b\a_s \rangle =-1$, 
thus $\b_1 +\b_2 -b\a_s \in A$ (since $A$ is closed), 
therefore $\b_1 + \b_2 \in \sigma_s \ast A$ (by Lemma 4.4).

\smallskip\noindent
{\it Case 3:} $\langle \a_s, \b_1 \rangle =0$ and $\langle 
\a_s, \b_2 \rangle =-b<0$. Then $\langle \a_s, \b_1+\b_2 
\rangle =-b<0$. Moreover, one has $\b_1, \b_2, \b_2+b\a_s 
\in A$ (by Lemma 4.4) and $\langle \b_1, \b_2 \rangle = 
\langle \b_1, \b_2+b\a_s \rangle =-1$, thus $\b_1+\b_2, \b_1 
+\b_2 +b\a_s \in A$ (since $A$ is closed), therefore 
$\b_1+\b_2 \in \sigma_s \ast A$ (by Lemma 4.4).

\smallskip\noindent
{\it Case 4:} $\langle \a_s, \b_1 \rangle =a>0$ and $\langle 
\a_s, \b_2 \rangle =b>0$. Then $\langle \a_s, \b_1+\b_2 
\rangle =a+b>0$. Moreover, one has $\b_1-a\a_s, \b_2-b\a_s 
\in A$ (by Lemma 4.4) and $\langle \b_1-a\a_s, \b_2-b\a_s 
\rangle =-1$, thus $\b_1 +\b_2 -(a+b) \a_s \in A$ (since $A$ 
is closed), therefore $\b_1+\b_2 \in \sigma_s \ast A$ (by 
Lemma 4.4).

\smallskip\noindent
{\it Case 5:} $\langle \a_s, \b_1 \rangle =a>0$ and $\langle 
\a_s, \b_2 \rangle =-b<0$. Note first that $\b_2, \b_2+b\a_s 
\in A$ (by Lemma 4.4), thus $\langle \b_2, \b_2+b\a_s 
\rangle = 2-b^2 \ge -1$ (since $A$ is closed), therefore 
$b=1$.

Suppose $a=1$. Then $\langle \a_s, \b_1+\b_2 \rangle =0$. 
One has $\b_1-\a_s, \b_2+\a_s \in A$ (by Lemma 4.4) and 
$\langle \b_1-\a_s, \b_2+\a_s \rangle =-1$, thus $\b_1+\b_2 
\in A$ (since $A$ is closed), therefore $\b_1+\b_2 \in 
\sigma_s \ast A$ (by Lemma 4.4).

Suppose $a\ge 2$. Then $\langle \a_s, \b_1+\b_2 \rangle =a-
1>0$. One has $\b_1-a\a_s, \b_2+\a_s \in A$ (by Lemma 4.4) 
and $\langle \b_1-a\a_s, \b_2+\a_s \rangle =-1$, thus 
$\b_1+\b_2 -(a-1)\a_s \in A$ (since $A$ is closed), 
therefore $\b_1+\b_2 \in \sigma_s \ast A$ (by Lemma 4.4).

\smallskip\noindent
{\it Case 6:} $\langle \a_s, \b_1 \rangle =-a<0$ and $\langle 
\a_s, \b_2 \rangle =-b<0$. Then $\b_1, \b_2+b\a_s \in A$ (by 
Lemma 4.4) and $\langle \b_1, \b_2+b\a_s \rangle = -1-ab <-
1$. This contradicts the definition of a closed subset, thus 
this case does not hold. \qed

\bigskip\noindent
{\bf Corollary 4.6.} {\it Let $A$ be a closed subset of 
$\Phi^+$ and $g \in G_\Gamma^+$. Then $g \ast A$ is also a 
closed subset.}

\bigskip\noindent
{\bf Lemma 4.7.} {\it Let $w \in W$ and $s \in S$ such that 
$l(sw) =l(w)-1$, and let $A$ be a closed subset of $\Phi^+$. 
One has $\Phi_w \subset \{ \a_s \} \cup s(A \setminus \{ 
\a_s \})$ if and only if $w< L( \sigma_s \cdot (\tau \circ 
C) (A))$.}

\bigskip\noindent
{\bf Proof.} The equality $l(sw)=l(w)-1$ implies, by 
Proposition 2.4, that $\Phi_w= \{ \a_s \} \sqcup s \cdot 
\Phi_{sw}$. So, the inclusion $\Phi_w \subset \{ \a_s \} 
\cup s(A \setminus \{ \a_s \})$ is equivalent to $s( 
\Phi_{sw}) \subset s(A \setminus \{ \a_s \})$, which is 
equivalent to $\Phi_{sw} \subset A$ (we may not have $\a_s 
\in \Phi_{sw}$ because $l(sw)<l(w)$). This inclusion is 
equivalent to $sw< C(A)$, which is equivalent to $\tau (sw) 
< (\tau \circ C)(A)$, which is equivalent to $\sigma_s \cdot 
\tau(sw) = \tau(w) < \sigma_s \cdot (\tau \circ C) (A)$, 
which is equivalent to $w< L( \sigma_s \cdot (\tau \circ C) 
(A))$. \qed

\bigskip\noindent
{\bf Lemma 4.8.} {\it Let $A$ and $B$ be two closed subsets 
of $\Phi^+$ and $s \in S$. If $\{ \a_s\} \subset B \subset 
\{ \a_s\} \cup s(A \setminus \{ \a_s \})$, then $B \subset 
\sigma_s \ast A$.}

\bigskip\noindent
{\bf Proof.} By Lemma 4.4, it suffices to show that: if $\b 
\in B$ is such that $\langle \a_s, \b \rangle =-a<0$, then 
$\b \in A$. One has $\a_s, \b \in B$ and $B$ is a closed 
subset, thus $\langle \a_s, \b \rangle =-1$ and $\b +\a_s \in B$. 
It follows that 
$s(\b) =\a_s+\b \in B \setminus \{ \a_s \} \subset s(A 
\setminus \{ \a_s \})$, thus $\b \in A$. \qed 

\bigskip\noindent
{\bf Lemma 4.9.} {\it Let $A$ be a closed subset of $\Phi^+$ 
and $g \in G_\Gamma^+$. Then
$$
C( g \ast A) = L( g \cdot (\tau \circ C) (A)).
$$}

\noindent
{\bf Proof.} We argue by induction on $l(g)$. Assume first 
that $l(g)=1$. Then $g=\sigma_s$ for some $s\in S$. Let $w_1 
= C( \sigma_s \ast A)$ and $w_2 = L( \sigma_s \cdot (\tau 
\circ C) (A))$. Since $\a_s \in \sigma_s \ast A$, one has 
$\Phi_s =\{ \a_s \} \subset \sigma_s \ast A$, thus $s < 
w_1$, namely, $l(sw_1) =l(w_1)-1$. Moreover, 
$\Phi_{w_1} \subset \sigma_s \ast A \subset \{ \a_s \} 
\cup s(A \setminus \{ \a_s \})$, thus, by Lemma~4.7, 
$w_1<w_2$. One has $\tau(s) = \sigma_s < \sigma_s \cdot 
(\tau \circ C)(A)$, thus $s<w_2$, namely, $l(sw_2)=l(w_2)-1$. 
By Lemma~4.7, it follows that $\{ \a_s \} \subset \Phi_{w_2} 
\subset \{ \a_s \} \cup s(A \setminus \{ \a_s \})$ and so, 
by Lemma 4.8, $\Phi_{w_2} \subset \sigma_s \ast A$. This 
implies that $w_2<w_1$.

Assume now that $l(g) \ge 2$. We write $g= \sigma_s g_1$ 
where $s \in S$ and $l(g_1) =l(g)-1$. Then, by induction and 
by Proposition 2.3,
$$
C(g \ast A) = C( \sigma_s \ast (g_1 \ast A)) = L( \sigma_s 
(\tau \circ C) (g_1 \ast A)) = L( \sigma_s (\tau \circ L) 
(g_1(\tau \circ C) (A))) = L(g (\tau \circ C) (A)). \quad 
\square
$$

\noindent
{\bf Definition.} Let $\CC$ denote the set of closed subsets 
of $\Phi^+$. For $w \in W$ we set
$$
\UU_w = \bigcup_{ A \in \CC,\ C(A) =w} U_A.
$$
Note that $\UU_w \neq \emptyset$ (since it contains 
$U_{\Phi_w}$), and one has $\UU_u \cap \UU_v = \emptyset$ if $u 
\neq v$.

\bigskip\noindent
{\bf Lemma 4.10.} {\it Let $g \in G_\Gamma^+$ and $w \in W$. 
Then
$$
\psi_\infty (g) \cdot \UU_w \subset \UU_{L(g \cdot \tau(w))}.
$$}

\noindent
{\bf Proof.} Let $A \in \CC$ such that $C(A)=w$. One has 
$\psi_\infty (g) \cdot U_A \subset U_{g \ast A}$, and, by Lemma 
4.9, $C(g \ast A) = L(g \cdot (\tau \circ C)(A)) = L(g \cdot 
\tau(w))$, thus $\psi_\infty (g) \cdot U_A \subset \UU_{L( g 
\cdot \tau(w))}$. This shows that $\psi_\infty (g) \cdot \UU_w 
\subset \UU_{L(g \cdot \tau(w))}$. \qed

\bigskip\noindent
{\bf Proof of Theorem 4.1.} Let $f,g \in G_\Gamma^+$ such 
that $\psi(f)= \psi(g)$. We write $f = \tau(u) f_1$ and $g = 
\tau(v) g_2$, where $u=L(f)$, $v=L(g)$, and $f_1, g_1 \in 
G_\Gamma^+$. Note that $u=1$ if and only if $f=1$, and $v=1$ 
if and only if $g=1$. Lemma 4.10 implies that $\psi_\infty(f) 
\cdot \UU_1 \subset \UU_{L(f)} = \UU_u$, and that $\psi_\infty (g) 
\cdot \UU_1 \subset \UU_{L(g)} = \UU_v$. Since $\psi_\infty (f) = 
\psi_\infty(g)$, and since $\UU_u \cap \UU_v = \emptyset$ if $u 
\neq v$, it follows that $u=v$.

We prove now that $f=g$ by induction on $l(f)$. If $l(f)=0$, 
then $f=1$, thus $u=v=1$, therefore $g=1$. Suppose $l(f)>0$. 
Then $l(f_1)<l(f)$ and
$$
\psi(f_1) = \psi(\tau(u))^{-1} \psi(f) = \psi(\tau(v))^{-1} 
\psi(g) = \psi(g_1).
$$
By the inductive hypothesis, it follows that $f_1=g_1$, thus 
$f =\tau(u)f_1 = \tau(v) g_1 =g$. \qed

\bigskip\bigskip\noindent
{\titre 5. The general case}

\bigskip\noindent
Now, we assume that $\Gamma$ is any Coxeter graph. The goal 
of this section is to prove the following.

\bigskip\noindent
{\bf Theorem 5.1.} {\it There exists an injective homomorphism $\phi: 
G_\Gamma^+ \to G_{\tilde \Gamma}^+$ from $G_\Gamma^+$ to an 
Artin monoid $G_{\tilde \Gamma}^+$
associated to a Coxeter graph $\tilde\Gamma$ 
of small type with no triangle.}

\bigskip
Since we already know by Corollary 4.2 that $\iota: 
G_{\tilde \Gamma}^+ \to G_{\tilde \Gamma}$ is injective, 
Theorem 5.1 finishes the proof of Theorem 1.1.

We start summarizing some well known properties of 
$G_\Gamma^+$ that can be found in [BS] and [Mic].

We say that $g \in G_\Gamma^+$ is a {\it common multiple} of a 
finite subset $F= \{ f_1, \dots, f_n \} \subset G_\Gamma^+$ 
if $f_i<g$ for all $i=1, \dots, n$. If $F= \{f_1, \dots, f_n 
\}$ has a common multiple, then it has a {\it least common 
multiple}, which is obviously unique, and which will be 
denoted by $f_1 \vee \dots \vee f_n$.

Let $s,t \in S$. The subset $\{ \sigma_s, \sigma_t \}$ has a 
common multiple if and only if $m_{s,t} < +\infty$. In that 
case, one has $\sigma_s \vee \sigma_t = \pprod( \sigma_s, 
\sigma_t; m_{s,t})$. More generally, for a subset $T \subset S$, the set 
$\Sigma_T =\{ \sigma_t; t\in T \}$ has a common multiple if 
and only if the subgroup $W_T$ of $W$ generated by $T$ is 
finite. In that case, the least common multiple of 
$\Sigma_T$ is denoted by $\Delta_T$. It is equal to 
$\tau(w_T)$, where $w_T$ denotes the element of maximal 
length in $W_T$. If $W$ is finite, namely, if $\Gamma$ is of 
spherical type, then we will denote by $\Delta= \Delta( 
\Gamma)$ the least common multiple of $\Sigma= \{ \sigma_s; 
s \in S\}$.

Let $T \subset S$ and $f,g \in G_\Gamma^+$. If $f$ and $g$ 
have a common multiple and both lie in the submonoid 
generated by $\Sigma_T = \{ \sigma_t; t \in T \}$, then $f 
\vee g$ also lies in this submonoid.

\bigskip\noindent
{\bf Definition.} Let $\Gamma$ and $\Gamma'$ be two Coxeter 
graphs, let $S$ be the set of vertices of $\Gamma$, and let 
$\phi: G_\Gamma^+ \to G_{\Gamma'}^+$ be a homomorphism. We 
say that $\phi$ {\it respects lcm's} if

\smallskip
$\bullet$ $\phi( \sigma_s) \neq 1$ for all $s \in S$;

\smallskip
$\bullet$ $\{ \phi( \sigma_s), \phi (\sigma_t) \}$ has a 
common multiple if and only if $m_{s,t}<+ \infty$;

\smallskip
$\bullet$ if $m_{s,t} < +\infty$, then $\phi( \sigma_s \vee 
\sigma_t) = \phi(\sigma_s) \vee \phi(\sigma_t)$.

\bigskip
The following theorem can be found in [Cri, Thm. 1.3].

\bigskip\noindent
{\bf Theorem 5.2} (Crisp). {\it If a homomorphism $\phi: 
G_\Gamma^+ \to G_{\Gamma'}^+$ between Artin monoids respects 
lcm's, then it is injective.}

\bigskip\noindent
{\bf Proof of Theorem 5.1.} Let $A_n$ be the Coxeter graph 
of Figure 1. Let $f,g$ be the elements of $G_{A_n}^+$ 
defined by $f= \sigma_1 \sigma_3 \sigma_5 \dots$ and $g = 
\sigma_2 \sigma_4 \sigma_6 \dots$. It is shown in [BS, Lem. 
5.8] that
$$
\pprod(f,g; n+1) = \pprod(g,f; n+1) = \Delta(A_n). \leqno(1)
$$

\bigskip
\epsfysize=0.5 truecm
\centerline{\epsfbox{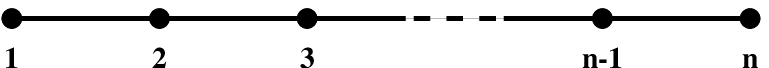}}
\bigskip
\centerline{Figure 1: The Coxeter graph $A_n$}

\bigskip
Let $m \ge 3$, and let $\Gamma(m)$ denote the Coxeter graph 
illustrated in Figure 2. It is a bipartite graph whose set 
of vertices is the disjoint union $I \sqcup J$, where $|I| = 
|J| = m-1$. As a Coxeter graph, $\Gamma(m)$ is the disjoint 
union of two copies of $A_{m-1}$. Let $f,g$ be the elements 
of $G_{\Gamma(m)}^+$ defined by $f= \prod_{i \in I} 
\sigma_i$ and $g = \prod_{j \in J} \sigma_j$. Then, by (1), 
one has
$$
\pprod (f,g;m) = \pprod(g,f;m) = \Delta( \Gamma(m)). 
\leqno(2)
$$

\bigskip
\epsfysize=3 truecm
\centerline{\epsfbox{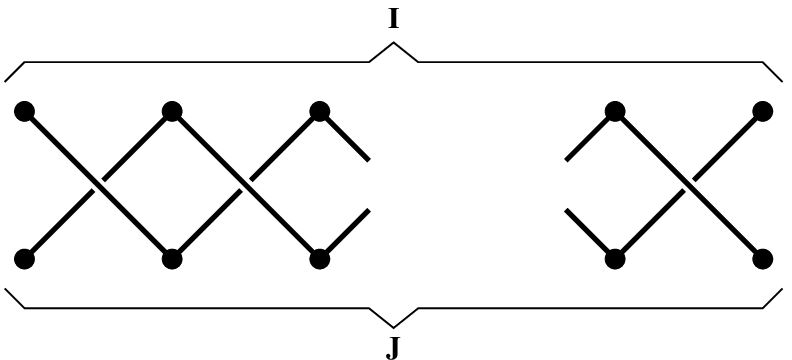}}
\bigskip
\centerline{Figure 2: The Coxeter graph $\Gamma(m)$}

\bigskip
Let $k \in \BN$. We denote by $k \Gamma(m)$ the 
disjoint union of $k$ copies of $\Gamma(m)$. It is a 
bipartite graph whose set of vertices is the disjoint union 
$kI \sqcup kJ$, where $kI$ denotes the disjoint union of $k$ 
copies of $I$, and $kJ$ denotes the disjoint union of $k$ 
copies of $J$. Let $f,g$ be the elements of $G_{k  
\Gamma(m)}^+$ defined by $f= \prod_{i \in kI} \sigma_i$ and 
$g= \prod_{j \in kJ} \sigma_j$. Then, by (2), one has
$$
\pprod (f,g;m) = \pprod(g,f;m) = \Delta (k \Gamma(m)). 
\leqno(3)
$$

Let $\Gamma( \infty)$ denote the Coxeter graph illustrated 
in Figure 3. It is bipartite graph whose set of vertices is 
the disjoint union $I \sqcup J$, where $I= \{ i_1, i_2 \}$ 
and $J= \{ j_1, j_2 \}$. Let $f,g$ be the elements of 
$G_{\Gamma( \infty)}^+$ defined by $f= \sigma_{i_1} 
\sigma_{i_2}$ and $g= \sigma_{j_1} \sigma_{j_2}$. A common 
multiple of $f$ and $g$ would be a common multiple of $\{ 
\sigma_{i_1}, \sigma_{i_2}, \sigma_{j_1}, \sigma_{j_2} \}$. 
But $\Gamma (\infty)$ is not of spherical type, thus such a 
common multiple may not exist. So, $f$ and $g$ have no 
common multiple.

Let $k \in \BN$. We denote by $k \Gamma (\infty)$ the 
disjoint union of $k$ copies of $\Gamma (\infty)$. It is a 
bipartite graph whose set of vertices is the disjoint union 
$kI \sqcup kJ$, where $kI$ denotes the disjoint union of $k$ 
copies of $I$, and $kJ$ denotes the disjoint union of $k$ 
copies of $J$. Let $f$ and $g$ be the elements of $G_{k 
\Gamma (\infty)}^+$ defined by $f= \prod_{i \in kI} 
\sigma_i$ and $g= \prod_{j \in kJ} \sigma_j$. Then, as 
before, $f$ and $g$ have no common multiple.

\vfill\eject

\bigskip
\epsfysize=2.5 truecm
\centerline{\epsfbox{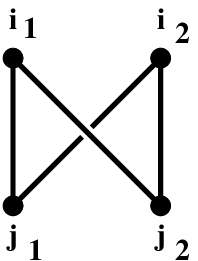}}
\bigskip
\centerline{Figure 3: The Coxeter graph $\Gamma(\infty)$}

\bigskip
Now, let $\Gamma$ be any Coxeter graph. Let $N$ be the least 
common multiple of $\{ m_{s,t}-1; s,t\in S, s\neq t, 
m_{s,t}< +\infty \}$. For all $s \in S$ we take an abstract 
set $I(s)$ with $2N$ elements. We construct a Coxeter graph 
of small type $\Gamma'$ as follows.

\smallskip
$\bullet$ The set of vertices of $\Gamma'$ is the disjoint 
union of the $I(s)$, $s \in S$.

\smallskip
$\bullet$ If $m_{s,t}=2$, then there is no edge joining two 
vertices in $I(s) \sqcup I(t)$.

\smallskip
$\bullet$ If $3 \le m_{s,t} < +\infty$, then the full 
subgraph of $\Gamma'$ generated by $I(s) \sqcup I(t)$ is 
isomorphic to 
\break
$\left( {2N \over m_{s,t}-1} \right) 
\Gamma (m_{s,t})$ with an isomorphism which takes $I(s)$ to 
$\left( {2N \over m_{s,t}-1} \right) I$ and $I(t)$ to 
$\left( {2N \over m_{s,t}-1} \right) J$.

\smallskip
$\bullet$ If $m_{s,t}= +\infty$, then the full subgraph of 
$\Gamma'$ generated by $I(s) \sqcup I(t)$ is isomorphic to 
$N \Gamma (\infty)$ with an isomorphism which takes 
$I(s)$ to $NI$ and $I(t)$ to $NJ$.

\smallskip
Such a graph always exists but is not unique in general. By 
the above considerations, there is a well defined 
homomorphism $\phi: G_\Gamma^+ \to G_{\Gamma'}^+$ which 
sends $\sigma_s$ to $\prod_{i \in I(s)} \sigma_i$ for all $s 
\in S$, and this homomorphism respects lcm's, so, is 
injective by Theorem 5.2. Note also that: if $\Gamma$ is of 
small type, then $\Gamma'$ is a bibartite graph, thus has no 
triangle.

So, applying twice the above construction, one gets a 
Coxeter graph $\tilde \Gamma$ of small type with no triangle 
and a monomorphism $\phi: G_\Gamma^+ \to G_{\tilde 
\Gamma}^+$. \qed

\bigskip\bigskip\noindent
{\titre 6. Two lemmas}

\bigskip\noindent
{\bf Lemma 3.3.} {\it Let $s \in S$ and $\b \in \Phi^+$ such that $\langle \a_s, \b \rangle =0$ 
and $\dpt (\b) \ge 2$. Then the definition of $T(s,\b)$ does not depend on the choice of the $t 
\in S$ such that $\dpt (t \cdot \b)= \dpt (\b)-1$.}

\bigskip\noindent
{\bf Proof.} We argue by induction on $\dpt (\b)$. We take $t,r \in S$, $t \neq r$, such that 
$\dpt( t \cdot \b) = \dpt (r \cdot \b) = \dpt (\b) -1$. By Proposition 2.6, we can write $ 
\langle \a_t, \b \rangle =b>0$ and $\langle \a_t, \b \rangle =c>0$.

\smallskip\noindent
{\it Case 1:} $\langle \a_s, \a_t \rangle =0$, $\langle \a_s, \a_r \rangle =0$, and $\langle 
\a_t, \a_r \rangle =0$.  We cannot have $\dpt (\b)=2$ in this case. Suppose $\dpt (\b) \ge 3$. 
Then, by induction, 
\bigskip
\centerline{\vbox{\halign{
\hfill$#$\ &$#$\hfill&\quad #\hfill\cr
y \cdot T(s, \b -b\a_t) &=\ y^2 \cdot T(s, \b -b\a_t -c\a_r) &by (D4)\cr
&=\ y \cdot T(s, \b -c\a_r) &by (D4)\cr}}}
\bigskip

\noindent
{\it Case 2:} $\langle \a_s, \a_t \rangle =0$, $\langle \a_s, \a_r \rangle =0$, and $\langle 
\a_t, \a_r \rangle =-1$. Suppose $\dpt(\b)=2$. Then $\b =\a_t +\a_r$, $b=c=1$, and
$$
y \cdot T(s, \b -\a_t) = y \cdot T(s, \a_r) =0=y \cdot T(s, \a_t) = y \cdot T(s, \b - \a_r) 
\quad {\rm by\ (D2)}
$$
We cannot have $\dpt (\b) =3$ in this case. Suppose $\dpt (\b) \ge 4$. Then, by induction
\bigskip
\centerline{\vbox{\halign{
\hfill$#$\ &$#$\hfill &\quad #\hfill\cr
y \cdot T(s, \b -b\a_t)&=\ y^2 \cdot T(s, \b -b\a_t -(b+c)\a_r) &by (D4)\cr
&=\ y^3 \cdot T(s, \b -(b+c)\a_t -(b+c)\a_r) &by (D4)\cr
&=\ y \cdot T(s, \b -c\a_r) &by symmetry\cr}}}
\bigskip

\noindent
{\it Case 3:} $\langle \a_s, \a_t \rangle =0$, $\langle \a_s, \a_r \rangle =-1$, and $\langle 
\a_t, \a_r \rangle =0$. We cannot have $\dpt(\b) \le 3$ in this case. Suppose $\dpt (\b) \ge 4$. 
Then, by induction,
\bigskip
\centerline{\vbox{\halign{
\hfill$#$ &$#$\hfill &\quad#\hfill\cr
y \cdot T(s, \b -b\a_t) &=\ y(y-1) \cdot T(s, \b -b\a_t -c\a_r) +y^2 \cdot T(r, \b -c\a_s -b\a_t 
-c\a_r) &by (D5)\cr
&=\ y^{\dpt (\b)-1} (y-1)^2 +y^2 \cdot T(r, \b -c\a_s -b\a_t -c\a_r) &by (D3)\cr
&=\ (y-1) \cdot T(s, \b -c\a_r) +y \cdot T(r, \b -c\a_s -c\a_r) &by (D3) and (D4)\cr}}}
\bigskip

\noindent
{\it Case 4:} $\langle \a_s, \a_t \rangle =0$, $\langle \a_s, \a_r \rangle =-1$, and $\langle 
\a_t, \a_r \rangle =-1$. We cannot have $\dpt (\b) \le 5$ in this case. Suppose $\dpt (\b) \ge 
6$. Then, by induction,
\bigskip
\centerline{\vbox{\halign{
\hfill$#$\ &$#$\hfill &\quad#\hfill\cr
&y \cdot T(s, \b -b\a_t)\cr
= &y(y-1) \cdot T(s, \b -b\a_t -(b+c)\a_r) +y^2 \cdot T(r, \b -(b+c)\a_s -b\a_t -(b+c)\a_r) &by 
(D5)\cr
= &y^{\dpt (\b)-1} (y-1)^2 +y^2(y-1) \cdot T(r, \b -(b+c)\a_s -(b+c)\a_t -(b+c)\a_r)\cr 
&\quad +y^3 \cdot T(t, \b -(b+c)\a_s -(b+c)\a_t -(b+2c)\a_r) &by (D3) and (D5)\cr
= &y^{\dpt(\b)-1} (y-1)^2 +y^{\dpt(\b)-2} (y-1)^2 +y^3 \cdot
T(t, \b -(b+c)\a_s -(b+c)\a_t -(b+2c)\a_r) &by (D3)\cr
= &y^{\dpt (\b)-1} (y-1)^2 +y(y-1) \cdot T(r, \b -c\a_s -(b+c)\a_t -c\a_r)\cr 
&\quad +y^2 \cdot T(t, \b -c\a_s -(b+c)\a_t -(b+2c)\a_r) &by (D3) and (D4)\cr
= &(y-1) \cdot T(s, \b -c\a_r) +y \cdot T(r, \b -c\a_s -c\a_r) &by (D3) and (D5)\cr}}}
\bigskip

\noindent
{\it Case 5:} $\langle \a_s, \a_t \rangle =-1$, $\langle \a_s, \a_r \rangle =-1$, and $\langle 
\a_t, \a_r \rangle =0$. We cannot have $\dpt (\b) \le 5$ in this case. Suppose $\dpt (\b) \ge 6$. 
Then, by induction,
\bigskip
\centerline{\vbox{\halign{
\hfill$#$\ &$#$\hfill &\quad #\hfill\cr
&(y-1) \cdot T(s, \b -b\a_t) +y \cdot T(t, \b -b\a_s -b\a_t)\cr
= &y^{\dpt(\b)-1} (y-1)^2 +y^2 \cdot T(t, \b -b\a_s -b\a_t -(b+c)\a_r) &by (D3) and (D4)\cr
= &y^{\dpt(\b)-1} (y-1)^2 +y^2(y-1) \cdot 
T(t, \b -(b+c)\a_s -b\a_t -(b+c)\a_r)\cr 
&\quad +y^3 \cdot T(s, \b -(b+c)\a_s -(b+c)\a_t -(b+c)\a_r) &by (D5)\cr
= &y^{\dpt(\b)-1} (y-1)^2 +y^{\dpt(\b)-2} (y-1)^2 +y^3 \cdot T(s, \b -(b+c)\a_s -(b+c)\a_t -
(b+c)\a_r) &by (D3)\cr
= &(y-1) \cdot T(s, \b -c\a_r) +y \cdot T(r, \b -c\a_s -c\a_r) &by symmetry.  \quad $\square$\cr}}}
\bigskip

\noindent
{\bf Lemma 3.4.} {\it Let $s \in S$ and $\b \in \Phi^+$ such that $\dpt (\b) \ge 2$ and $\langle 
\a_s, \b \rangle =-a<0$. Then the definition of $T(s,\b)$ does not depend on the choice of the $t 
\in S$ such that $\dpt (t \cdot \b) = \dpt(\b)-1$.}

\bigskip\noindent
{\bf Proof.} We argue by induction on $\dpt (\b)$. We take $t,r \in S$, $t \neq r$, such that 
$\dpt (t \cdot \b) = \dpt (r \cdot \b) = \dpt(\b)-1$. By Proposition 2.6, we can write $\langle 
\a_t, \b \rangle =b>0$ and $\langle \a_r, \b \rangle =c>0$.

\smallskip\noindent
{\it Case 1:} $\langle \a_s, \a_t \rangle =0$, $\langle \a_s, \a_r \rangle =0$, and $\langle 
\a_t, \a_r \rangle =0$. We cannot have $\dpt (\b)=2$ in this case. Suppose $\dpt(\b) \ge 3$. 
Then, by induction,
\bigskip
\centerline{\vbox{\halign{
\hfill$#$\ &$#$\hfill&\quad #\hfill\cr
y \cdot T(s, \b -b\a_t) &=\ y^2 \cdot T(s, \b -b\a_t -c\a_r) &by (D6)\cr
&=\ y \cdot T(s, \b -c\a_r) &by (D6)\cr}}}
\bigskip

\noindent
{\it Case 2:} $\langle \a_s, \a_t \rangle =0$, $\langle \a_s, \a_r \rangle =0$, and $\langle 
\a_t, \a_r \rangle =-1$. We cannot have $\dpt(\b) \le 3$ in this case. Suppose $\dpt (\b) \ge 4$. 
Then, by induction,
\bigskip
\centerline{\vbox{\halign{
\hfill$#$\ &$#$\hfill&\quad #\hfill\cr
y \cdot T(s, \b -b\a_t) &=\ y^2 \cdot T(s, \b -b\a_t -(b+c)\a_r) &by (D6)\cr
&=\ y^3 \cdot T(s, \b -(b+c)\a_t -(b+c)\a_r) &by (D6)\cr
&=\ y \cdot T(s, \b -c\a_r) &by symmetry\cr}}}
\bigskip

\noindent
{\it Case 3:} $\langle \a_s, \a_t \rangle =0$, $\langle \a_s, \a_r \rangle =-1$, $\langle \a_t, 
\a_r \rangle =0$, and $c>a$. We cannot have $\dpt(\b) \le 3$ in this case. Suppose $\dpt (\b) \ge 
4$. Then, by induction,
\bigskip
\centerline{\vbox{\halign{
\hfill$#$\ &$#$\hfill&\quad #\hfill\cr
&y \cdot T(s, \b -b\a_t)\cr
= &y(y-1) \cdot T(s, \b -b\a_t -c\a_r) +y^2 \cdot T(r, \b -(c-a)\a_s -b\a_t -c\a_r) &by (D7)\cr 
= &y^{\dpt(\b)-1} (y-1)^2 +y^2 \cdot T(r, \b -(c-a)\a_s -b\a_t -c\a_r) &by (D3)\cr 
= &(y-1) \cdot T(s, \b -c\a_r) +y \cdot T(r, \b -(c-a)\a_s -c\a_r) &by (D3) and (D6)\cr}}}
\bigskip

\noindent
{\it Case 4:} $\langle \a_s, \a_t \rangle =0$, $\langle \a_s, \a_r \rangle =-1$, $\langle \a_t, 
\a_r \rangle =0$, and $c=a$. We cannot have $\dpt(\b) =2$ in this case. Suppose $\dpt (\b) \ge 
3$. Then, by induction,
\bigskip
\centerline{\vbox{\halign{
\hfill$#$\ &$#$\hfill&\quad #\hfill\cr
&y \cdot T(s, \b -b\a_t)\cr
= &y \cdot T(r, \b -b\a_t -a\a_r) +y(y-1) \cdot T(s, \b -b\a_t -a\a_r) &by (D8)\cr
= &T(r, \b -a\a_r) +(y-1) \cdot T(s, \b -a\a_r) &by (D6) and (D4)\cr}}}
\bigskip

\noindent
{\it Case 5:} $\langle \a_s, \a_t \rangle =0$, $\langle \a_s, \a_r \rangle =-1$, $\langle \a_t, 
\a_r \rangle =0$, and $c<a$. We cannot have $\dpt(\b) =2$ in this case. Suppose $\dpt (\b) \ge 
3$. Then, by induction,
\bigskip
\centerline{\vbox{\halign{
\hfill$#$\ &$#$\hfill&\quad #\hfill\cr
&y \cdot T(s, \b -b\a_t)\cr
= &y^2 \cdot T(s, \b -b\a_t -c\a_r) +y \cdot T(r, \b -b\a_t -c\a_r) +y^{\dpt(\b)-1} (1-y) &by 
(D9)\cr
= &y \cdot T(s, \b -c\a_r) +T(r, \b -c\a_r) +y^{\dpt(\b)-1} (1-y) &by (D6)\cr}}}
\bigskip

\noindent
{\it Case 6:} $\langle \a_s, \a_t \rangle =0$, $\langle \a_s, \a_r \rangle =-1$, $\langle \a_t, 
\a_r \rangle =-1$, and $c>a$. We cannot have $\dpt (\b) \le 5$ in this case. Suppose $\dpt (\b) 
\ge 6$. Then, by induction,
\bigskip
\centerline{\vbox{\halign{
\hfill$#$\ &$#$\hfill&\quad #\hfill\cr
&y \cdot T(s, \b -b\a_t)\cr
= &y(y-1) \cdot T(s, \b -b \a_t -(b+c)\a_r) +y^2 \cdot T(r, \b -(b+c-a)\a_s -b\a_t -(b+c)\a_r) 
&by (D7)\cr
= &y^{\dpt(\b)-1} (y-1)^2 +y^2(y-1) \cdot T(r, \b -(b+c-a)\a_s -(b+c)\a_t -(b+c)\a_r)\cr 
&\quad
+y^3 \cdot T(t, \b -(b+c-a)\a_s -(b+c)\a_t -(b+2c-a)\a_r) &by (D3) and (D7)\cr
= &y^{\dpt(\b)-1} (y-1)^2 + y^{\dpt(\b)-2} (y-1)^2\cr 
&\quad
+y^3 \cdot T(t, \b -(b+c-a)\a_s -(b+c)\a_t -(b+2c-a)\a_r) &by (D3)\cr
= &y^{\dpt(\b)-1} (y-1)^2 + y(y-1) \cdot T(r, \b -(c-a)\a_s -(b+c)\a_t -c\a_r)\cr 
&\quad
+y^2 \cdot T(t, \b -(c-a)\a_s -(b+c)\a_t -(b+2c-a)\a_r) &by (D3) and (D6)\cr
= & (y-1) \cdot T(s, \b -c\a_r) +y \cdot T(r, \b -(c-a)\a_s  -c\a_r) &by (D3) and (D7)\cr}}}
\bigskip

\noindent
{\it Case 7:} $\langle \a_s, \a_t \rangle =0$, $\langle \a_s, \a_r \rangle =-1$, $\langle \a_t, 
\a_r \rangle =-1$, and $c=a$. Suppose $\dpt(\b)=2$. Then $a=b=c=1$, $\b =\a_t +\a_r$, and
\bigskip
\centerline{\vbox{\halign{
\hfill$#$&$#$\hfill&\quad#\hfill\cr
y \cdot T(s, \b -\a_t) =&\ y \cdot T(s, \a_r) =0 &by (D2)\cr
T(r, \b -\a_r) +(y-1) \cdot T(s, \b -\a_r) =&\ T(r, \a_t) + (y-1) \cdot T(s, \a_t) =0 &by 
(D2)\cr}}}
\bigskip\noindent
We cannot have $\dpt (\b) \in \{ 3,4 \}$ in this case. 
Suppose $\dpt(\b) \ge 5$. Then, by induction,
\bigskip
\centerline{\vbox{\halign{
\hfill$#$\ &$#$\hfill&\quad#\hfill\cr
&y \cdot T(s, \b -b\a_t)\cr
= &y(y-1) \cdot T(s, \b -b\a_t -(b+a)\a_r) +y^2 \cdot T(r, \b -b\a_s -b\a_t -(b+a)\a_r) &by 
(D7)\cr
= &y^{\dpt(\b)-1} (y-1)^2 +y^2 \cdot T(t, \b -b\a_s -(b+a)\a_t -(b+a)\a_r)\cr 
&\quad
+y^2(y-1) \cdot T(r, \b -b\a_s -(b+a)\a_t -(b+a)\a_r) &by (D3) and (D8)\cr
= &y^{\dpt(\b)-2} (y-1)^2 +y^2 \cdot T(t, \b -b\a_s -(b+a)\a_t -(b+a)\a_r) + y^{\dpt(\b)-2} (y-
1)^3\cr 
&\quad
+y^2(y-1) \cdot T(r, \b -b\a_s -(b+a)\a_t -(b+a)\a_r)\cr
= &y^{\dpt(\b)-2} (y-1)^2 +y^2 \cdot T(t, \b -b\a_s -(b+a)\a_t -(b+a)\a_r)\cr 
&\quad
+ y(y-1)^2 \cdot T(s, \b -(b+a)\a_t -(b+a)\a_r)\cr
&\quad
+y^2(y-1) \cdot T(r, \b -b\a_s -(b+a)\a_t -(b+a)\a_r) &by (D3)\cr
= &y^{\dpt(\b)-2} (y-1)^2 +y^2 \cdot T(t, \b -b\a_s -(b+a)\a_t -(b+a)\a_r)\cr 
&\quad
+y(y-1) \cdot T(s, \b -(b+a)\a_t -a\a_r) &by (D5)\cr
= &(y-1) \cdot T(r, \b -(b+a)\a_t -a\a_r) +y \cdot T(t, \b -(b+a)\a_t -(b+a)\a_r)\cr 
&\quad
+ y(y-1) \cdot T(s, \b -(b+a)\a_t -a\a_r) &by (D3) and (D6)\cr
= &T(r, \b -a\a_r)+(y-1) \cdot T(s, \b -a\a_r) &by (D7) and (D4)\cr}}}
\bigskip

\noindent
{\it Case 8:} $\langle \a_s, \a_t \rangle =0$, $\langle \a_s, \a_r \rangle =-1$, $\langle \a_t, 
\a_r \rangle =-1$, $c<a$, and $b+c>a$. We cannot have $\dpt (\b) \le 4$ in this case. Suppose 
$\dpt (\b) \ge 5$. Then, by induction,
\bigskip
\centerline{\vbox{\halign{
\hfill$#$\ &$#$\hfill&\quad#\hfill\cr
&y \cdot T(s, \b -b\a_t)\cr
= &y(y-1) \cdot T(s, \b -b\a_t -(b+c)\a_r) +y^2 \cdot T(r, \b -(b+c-a)\a_s -b\a_t -(b+c)\a_r) &by 
(D7)\cr
= &y^{\dpt(\b)-1} (y-1)^2 +y^3 \cdot T(r, \b -(b+c-a)\a_s -(b+c)\a_t -(b+c)\a_r)\cr 
&\quad
+y^2 \cdot T(t, \b -(b+c-a)\a_s -(b+c)\a_t -(b+c)\a_r) +y^{\dpt(\b)-2} (1-y) &by (D3) and (D9)\cr
= &y^{\dpt(\b)-1} (y-1)^2 +y^3 \cdot T(r, \b -(b+c-a)\a_s -(b+c)\a_t -(b+c)\a_r)\cr 
&\quad
+y^{\dpt(\b)-2} (y-1)^2\cr 
&\quad
+y^2 \cdot T(t, \b -(b+c-a)\a_s -(b+c)\a_t -(b+c)\a_r) +y^{\dpt(\b)-1} (1-y)\cr
= &y^2(y-1) \cdot T(s, \b -(b+c)\a_t -(b+c)\a_r)\cr 
&\quad
+y^3 \cdot T(r, \b -(b+c-a)\a_s -(b+c)\a_t -(b+c)\a_r) +y^{\dpt(\b)-2} (y-1)^2\cr 
&\quad
+y \cdot T(t, \b -(b+c)\a_t -(b+c)\a_r) +y^{\dpt(\b)-1} (1-y) &by (D3) and (D6)\cr
= &y^2 \cdot T(s, \b -(b+c)\a_t -c\a_r) +(y-1) \cdot T(r, \b -(b+c)\a_t -c\a_r)\cr 
&\quad
+y \cdot T(t, \b -(b+c)\a_t -(b+c)\a_r) +y^{\dpt(\b)-1} (1-y) &by (D7) and (D3)\cr
= &y \cdot T(s, \b -c\a_r) +T(r, \b -c\a_r) +y^{\dpt(\b)-1} (1-y) &by (D6) and (D7)\cr}}}
\bigskip

\noindent
{\it Case 9:} $\langle \a_s, \a_t \rangle =0$, $\langle \a_s, \a_r \rangle =-1$, $\langle \a_t, 
\a_r \rangle =-1$, $c<a$, and $b+c=a$. We cannot have $\dpt(\b) \le 3$ in this case. Suppose 
$\dpt(\b) \ge 4$. Then, by induction,
\bigskip
\centerline{\vbox{\halign{
\hfill$#$\ &$#$\hfill&\quad#\hfill\cr
&y \cdot T(s, \b -b\a_t)\cr
= &y \cdot T(r, \b -b\a_t -a\a_r) +y(y-1) \cdot T(s, \b -b\a_t -a\a_r) &by (D8)\cr
= &y^2 \cdot T(r, \b -a\a_t -a\a_r) +y \cdot T(t, \b -a\a_t -a\a_r) +y^{\dpt(\b)-2} (1-y)\cr 
&\quad
+y^2(y-1) \cdot T(s, \b -a\a_t -a\a_r) &by (D9) and (D4)\cr
= &y^2 \cdot T(r, \b -a\a_t -a\a_r) +y^2(y-1) \cdot T(s, \b -a\a_t -a\a_r) +y^{\dpt(\b)-2} (y-
1)^2\cr 
&\quad
+y \cdot T(t, \b -a\a_t -a\a_r) +y^{\dpt(\b)-1} (1-y)\cr
= &y^2 \cdot T(s, \b -a\a_t -c\a_r) +(y-1) \cdot T(r, \b -a\a_t -c\a_r)\cr 
&\quad
+y \cdot T(t, \b -a\a_t -a\a_r) +y^{\dpt(\b)-1} (1-y) &by (D8) and (D3)\cr
= &y \cdot T(s, \b -c\a_r) +T(r, \b -c\a_r) +y^{\dpt(\b)-1} (1-y) &by (D6) and (D7)\cr}}}
\bigskip

\noindent
{\it Case 10:} $\langle \a_s, \a_t \rangle =0$, $\langle \a_s, \a_r \rangle =-1$, $\langle \a_t, 
\a_r \rangle =-1$, $c<a$, and $b+c<a$. We cannot have $\dpt(\b) \le 3$ in this case. Suppose 
$\dpt(\b) \ge 4$. Then, by induction,
\bigskip
\centerline{\vbox{\halign{
\hfill$#$\ &$#$\hfill&\quad#\hfill\cr
&y \cdot T(s, \b -b\a_t)\cr
= &y^2 \cdot T(s, \b -b\a_t -(b+c)\a_r) +y \cdot T(r, \b -b\a_t -(b+c)\a_r) +y^{\dpt(\b)-1} (1-y) &by (D9)\cr
= &y^3 \cdot T(s, \b -(b+c)\a_t -(b+c)\a_r) +y^2 \cdot T(r, \b -(b+c)\a_t -(b+c)\a_r)\cr 
&\quad
+y \cdot T(t, \b -(b+c)\a_t -(b+c)\a_r) +y^{\dpt(\b)-2} (1-y) +y^{\dpt(\b)-1} (1-y) &by (D6) and 
(D9)\cr
= &y^3 \cdot T(s, \b -(b+c)\a_t -(b+c)\a_r) +y^2 \cdot T(r, \b -(b+c)\a_t -(b+c)\a_r)\cr 
&\quad
+y^{\dpt(\b)-1} (1-y) +y^{\dpt(\b)-2} (y-1)^2\cr 
&\quad
+y \cdot T(t, \b -(b+c)\a_t -(b+c)\a_r) +y^{\dpt(\b)-1} (1-y)\cr
= &y^2 \cdot T(s, \b -(b+c)\a_t -c\a_r) +(y-1) \cdot T(r, \b -(b+c)\a_t -c\a_r)\cr 
&\quad
+y \cdot T(t, \b -(b+c)\a_t -(b+c)\a_r) +y^{\dpt(\b)-1} (1-y) &by (D9) and (D3)\cr
= &y \cdot T(s, \b -c\a_r) +T(r, \b -c\a_r) +y^{\dpt(\b)-1} (1-y) &by (D6) and (D7)\cr}}}
\bigskip

\noindent
{\it Case 11:} $\langle \a_s, \a_t \rangle =-1$, $\langle \a_s, \a_r \rangle =-1$, $\langle \a_t, 
\a_r \rangle =0$, $b>a$, and $c>a$. We cannot have $\dpt(\b) \le 5$ in this case. Suppose 
$\dpt(\b) \ge 6$. Then, by induction,
\bigskip
\centerline{\vbox{\halign{
\hfill$#$\ &$#$\hfill&\quad#\hfill\cr
&(y-1) \cdot T(s, \b -b\a_t) +y \cdot T(t, \b -(b-a)\a_s -b\a_t)\cr
= &y^{\dpt(\b)-1} (y-1)^2 +y^2 \cdot T(t, \b -(b-a)\a_s -b\a_t -(b+c-a)\a_r) &by (D3) and (D6)\cr
= &y^{\dpt(\b)-1} (y-1)^2 +y^2(y-1) \cdot T(t, \b -(b+c-a)\a_s -b\a_t -(b+c-a)\a_r)\cr  
&\quad
+y^3 \cdot T(s, \b -(b+c-a)\a_s -(b+c-a)\a_t -(b+c-a)\a_r) &by (D7)\cr
= &y^{\dpt(\b)-1} (y-1)^2 +y^{\dpt(\b)-2} (y-1)^2\cr 
&\quad
+y^3 \cdot T(s, \b -(b+c-a)\a_s -(b+c-a)\a_t -(b+c-a)\a_r) &by (D3)\cr
= &(y-1) \cdot T(s, \b -c\a_r) +y \cdot T(r, \b -(c-a)\a_s -c\a_r) &by symmetry\cr}}}
\bigskip

\noindent
{\it Case 12:} $\langle \a_s, \a_t \rangle =-1$, $\langle \a_s, \a_r \rangle =-1$, $\langle \a_t, 
\a_r \rangle =0$, and $b>c=a$. We cannot have $\dpt(\b) \le 4$ in this case. Suppose $\dpt(\b) 
\ge 5$. Then, by induction,
\bigskip
\centerline{\vbox{\halign{
\hfill$#$\ &$#$\hfill&\quad#\hfill\cr
&(y-1) \cdot T(s, \b -b\a_t) +y \cdot T(t, \b -(b-a)\a_s -b\a_t)\cr
= &y^{\dpt(\b)-1} (y-1)^2 +y^2 \cdot T(t, \b -(b-a)\a_s -b\a_t-b\a_r) &by (D3) and (D6)\cr
= &y^{\dpt(\b)-1} (y-1)^2 +y^2 \cdot T(s, \b -b\a_s -b\a_t-b\a_r) +y^2(y-1) \cdot T(t, \b -b\a_s 
-b\a_t-b\a_r) &by (D8)\cr
= &y^{\dpt(\b)-2} (y-1)^2 +y^2 \cdot T(s, \b -b\a_s -b\a_t-b\a_r) +y^{\dpt(\b)-2} (y-1)^3\cr 
&\quad
+y^2(y-1) \cdot T(t, \b -b\a_s -b\a_t-b\a_r)\cr
= &y(y-1) \cdot T(r, \b -b\a_s -b\a_t-a\a_r) +y^2 \cdot T(s, \b -b\a_s -b\a_t-b\a_r) 
+y^{\dpt(\b)-2} (y-1)^3\cr 
&\quad
+y(y-1) \cdot T(t, \b -b\a_s -b\a_t-a\a_r) &by (D3) and (D4)\cr
= &y \cdot T(r, \b -b\a_t-a\a_r) +(y-1)^2 \cdot T(s, \b -b\a_t-a\a_r)\cr 
&\quad
+y(y-1) \cdot T(t, \b -b\a_s -b\a_t-a\a_r) &by (D7) and (D3)\cr
= &T(r, \b -a\a_r) +(y-1) \cdot T(s, \b -a\a_r) &by (D6) and (D5)\cr}}}
\bigskip

\noindent
{\it Case 13:} $\langle \a_s, \a_t \rangle =-1$, $\langle 
\a_s, \a_r \rangle =-1$, $\langle \a_t, \a_r \rangle =0$, and 
$b>a>c$. We cannot have $\dpt(\b) \le 3$ in this case. 
Suppose $\dpt(\b) \ge 4$. Then, by induction,
\bigskip
\centerline{\vbox{\halign{
\hfill$#$\ &$#$\hfill&\quad#\hfill\cr
&(y-1) \cdot T(s, \b -b\a_t) +y \cdot T(t, \b -(b-a)\a_s -
b\a_t)\cr
= &y^{\dpt(\b)-1} (y-1)^2 +y^2 \cdot T(t, \b -(b-a)\a_s -
b\a_t -(b+c-a)\a_r) &by (D3) and (D6)\cr
= &y^{\dpt(\b)-1} (y-1)^2 +y^3 \cdot T(t, \b -(b+c-a)\a_s -
b\a_t -(b+c-a)\a_r)\cr 
&\quad
+y^2 \cdot T(s, \b -(b+c-a)\a_s -b\a_t -(b+c-a)\a_r) 
+y^{\dpt(\b)-2} (1-y) &by (D9)\cr
= &y^{\dpt(\b)-1} (y-1)^2 +y^3 \cdot T(t, \b -(b+c-a)\a_s -
b\a_t -(b+c-a)\a_r) +y^{\dpt(\b)-2} (y-1)^2\cr 
&\quad
+y^2 \cdot T(s, \b -(b+c-a)\a_s -b\a_t -(b+c-a)\a_r) 
+y^{\dpt(\b)-1} (1-y)\cr
= &y^{\dpt(\b)-1} (y-1)^2 +y^2 \cdot T(t, \b -(b+c-a)\a_s -
b\a_t -c\a_r)\cr 
&\quad
+y(y-1) \cdot T(r, \b -(b+c-a)\a_s -b\a_t -c\a_r)\cr 
&\quad
+y^2 \cdot T(s, \b -(b+c-a)\a_s -b\a_t -(b+c-a)\a_r) 
+y^{\dpt(\b)-1} (1-y) &by (D6) and (D3)\cr
= &y(y-1) \cdot T(s, \b -b\a_t -c\a_r) +y^2 \cdot T(t, \b -
(b+c-a)\a_s -b\a_t -c\a_r)\cr 
&\quad
+y \cdot T(r, \b -b\a_t -c\a_r) +y^{\dpt(\b)-1} (1-y) &by 
(D3) and (D7)\cr
= &y \cdot T(s, \b -c\a_r) + T(r, \b -c\a_r) +y^{\dpt(\b)-1} 
(1-y) &by (D7) and (D6)\cr}}}
\bigskip

\noindent
{\it Case 14:} $\langle \a_s, \a_t \rangle =-1$, $\langle 
\a_s, \a_r \rangle =-1$, $\langle \a_t, \a_r \rangle =0$, and 
$a=b=c$. We cannot have $\dpt(\b) \le 3$ in this case. 
Suppose $\dpt(\b) \ge 4$. Then, by induction,
\bigskip
\centerline{\vbox{\halign{
\hfill$#$\ &$#$\hfill&\quad#\hfill\cr
&T(t, \b -a\a_t) +(y-1) \cdot T(s, \b -a\a_t)\cr
= &y \cdot T(t, \b -a\a_t -a\a_r) +(y-1)^2 \cdot T(s, \b -
a\a_t-a\a_r)\cr 
&\quad
+y(y-1) \cdot T(r, \b -a\a_s -a\a_t-a\a_r) &by (D6) and 
(D5)\cr
= &y \cdot T(s, \b -a\a_s -a\a_t-a\a_r) +y(y-1) \cdot T(t, \b 
-a\a_s -a\a_t -a\a_r)\cr 
&\quad
+y^{\dpt(\b)-2} (y-1)^3 +y(y-1) \cdot T(r, \b -a\a_s -a\a_t-
a\a_r) &by (D8) and (D3)\cr
= &T(r, \b -a\a_r) +(y-1) \cdot T(s, \b -a\a_r) &by 
symmetry\cr}}}
\bigskip

\noindent
{\it Case 15:} $\langle \a_s, \a_t \rangle =-1$, $\langle 
\a_s, \a_r \rangle =-1$, $\langle \a_t, \a_r \rangle =0$, and 
$a=b>c$. We cannot have $\dpt(\b) \le 3$ in this case. 
Suppose $\dpt(\b) \ge 4$. Then, by induction,
\bigskip
\centerline{\vbox{\halign{
\hfill$#$\ &$#$\hfill&\quad#\hfill\cr
&T(t, \b -a\a_t) +(y-1) \cdot T(s, \b -a\a_t)\cr
= &y \cdot T(t, \b -a\a_t -c\a_r) +(y-1)^2 \cdot T(s, \b -
a\a_t -c\a_r)\cr 
&\quad
+y(y-1) \cdot T(r, \b -c\a_s -a\a_t -c\a_r) &by (D6) and 
(D5)\cr
= &y^2 \cdot T(t, \b -c\a_s -a\a_t -c\a_r) +y \cdot T(s, \b -
c\a_s -a\a_t -c\a_r) +y^{\dpt(\b)-2} (1-y)\cr 
&\quad
+y^{\dpt(\b)-2} (y-1)^3 +y(y-1) \cdot T(r, \b -c\a_s -a\a_t -
c\a_r) &by (D9) and (D3)\cr
= &y^{\dpt(\b)-1} (y-1)^2 +y^2 \cdot T(t, \b -c\a_s -a\a_t -
c\a_r) +y \cdot T(s, \b -c\a_s -a\a_t -c\a_r)\cr 
&\quad
+y(y-1) \cdot T(r, \b -c\a_s -a\a_t -c\a_r) +y^{\dpt(\b)-1} 
(1-y)\cr
= &y(y-1) \cdot T(s, \b -a\a_t -c\a_r) +y^2 \cdot T(t, 
\b -c\a_s -a\a_t -c\a_r)\cr 
&\quad
+y \cdot T(r, \b -a\a_t -c\a_r) +y^{\dpt(\b)-1} (1-y) &by 
(D3) and (D8)\cr
= &y \cdot T(s, \b -c\a_r) +T(r, \b -c\a_r) +y^{\dpt(\b)-1} 
(1-y) &by (D7) and (D6)\cr}}}
\bigskip

\noindent
{\it Case 16:} $\langle \a_s, \a_t \rangle =-1$, $\langle 
\a_s, \a_r \rangle =-1$, $\langle \a_t, \a_r \rangle =0$, 
$a>b$, $a>c$, and $b+c>a$. We cannot have $\dpt(\b) \le 3$ in 
this case. Suppose $\dpt(\b) \ge 4$. Then, by induction,
\bigskip
\centerline{\vbox{\halign{
\hfill$#$\ &$#$\hfill&\quad#\hfill\cr
&y \cdot T(s, \b -b\a_t) +T(t, \b -b\a_t) +y^{\dpt(\b)-1} (1-
y)\cr
= &y(y-1) \cdot T(s, \b -b\a_t-c\a_r) + y^2 \cdot T(r, \b-
(b+c-a)\a_s -b\a_t-c\a_r)\cr
&\quad
+ y \cdot T(t, \b -b\a_t -c\a_r) +y^{\dpt(\b)-1} (1-y) &by 
(D7) and (D6)\cr
= &y^{\dpt(\b)-1} (y-1)^2 + y^2 \cdot T(r, \b-(b+c-a)\a_s -
b\a_t-c\a_r)\cr
&\quad
+ y^2 \cdot T(t, \b -(b+c-a)\a_s -b\a_t -c\a_r)\cr 
&\quad
+y \cdot T(s, \b -(b+c-a)\a_s -b\a_t -c\a_r) +y^{\dpt(\b)-2} 
(1-y) +y^{\dpt(\b)-1} (1-y) &by (D3) and (D9)\cr
 = &y \cdot T(s, \b -c\a_r) +T(r, \b -c\a_r) +y^{\dpt(\b)-1} 
(1-y) &by symmetry\cr}}}
\bigskip

\noindent
{\it Case 17:} $\langle \a_s, \a_t \rangle =-1$, $\langle 
\a_s, \a_r \rangle =-1$, $\langle \a_t, \a_r \rangle =0$, and 
$a=b+c$. We cannot have $\dpt(\b) =2$ in this case. Suppose 
$\dpt(\b) \ge 3$. Then, by induction,
\bigskip
\centerline{\vbox{\halign{
\hfill$#$\ &$#$\hfill&\quad#\hfill\cr
&y \cdot T(s, \b -b\a_t) +T(t, \b -b\a_t) +y^{\dpt(\b)-1} (1-
y)\cr
= &y \cdot T(r, \b -b\a_t -c\a_r) +y(y-1) \cdot T(s, \b -
b\a_t -c\a_r)\cr
&\quad
+y \cdot T(t, \b -b\a_t -c\a_r) +y^{\dpt(\b)-1} (1-y) &by 
(D8) and (D6)\cr
 = &y \cdot T(s, \b -c\a_r) +T(r, \b -c\a_r) +y^{\dpt(\b)-1} 
(1-y) &by symmetry\cr}}}
\bigskip

\noindent
{\it Case 18:} $\langle \a_s, \a_t \rangle =-1$, $\langle 
\a_s, \a_r \rangle =-1$, $\langle \a_t, \a_r \rangle =0$, and 
$a>b+c$. We cannot have $\dpt(\b) =2$ in this case. Suppose 
$\dpt(\b) \ge 3$. Then, by induction,
\bigskip
\centerline{\vbox{\halign{
\hfill$#$\ &$#$\hfill&\quad#\hfill\cr
&y \cdot T(s, \b -b\a_t) +T(t, \b -b\a_t) +y^{\dpt(\b)-1} (1-
y)\cr
= &y^2 \cdot T(s, \b -b\a_t -c\a_r)+ y \cdot T(r, \b -b\a_t -
c\a_r) +y^{\dpt(\b)-1} (1-y)\cr 
&\quad
+y \cdot T(t, \b -b\a_t -c\a_r) +y^{\dpt(\b)-1} (1-y) &by 
(D9) and (D6)\cr
= &y \cdot T(s, \b -c\a_r) +T(r, \b -c\a_r) +y^{\dpt(\b)-1} 
(1-y) &by symmetry. \quad $\square$\cr}}}

\bigskip\bigskip\noindent
{\titre References}

\bigskip
\item{[Alt]}
J. Altobelli,
{\it The word problem for Artin groups of FC type},
J. Pure Appl. Algebra {\bf 129} (1998), 1--22.

\smallskip
\item{[Big]}
S. Bigelow,
{\it Braid groups are linear},
J. Amer. Math. Soc., to appear.

\smallskip
\item{[BS]}
E. Brieskorn, K. Saito,
{\it Artin-Gruppen und Coxeter-Gruppen},
Invent. Math. {\bf 17} (1972), 245--271.

\smallskip
\item{[BH]}
B. Brink, R.B. Howlett,
{\it A finiteness property and an automatic structure 
for Coxeter groups},
Math. Ann. {\bf 296} (1993), 179--190.

\smallskip
\item{[Cha]}
R. Charney,
{\it Injectivity of the positive monoid for some infinite
type Artin groups},
Cossey, John (ed.) et al., Geometric group theory down under.
Proceedings of a special year in geometric
group theory, Camberra, Australia, July 14--19, 1996.
De Gruyter, Berlin, 1999, pp. 103--118.

\smallskip
\item{[ChP]}
J.R. Cho, S.J. Pride,
{\it Embedding semigroups into groups, and the asphericity
of semigroups},
Int. J. Algebra Comput. {\bf 3} (1993), 1--13.

\smallskip
\item{[CW]}
A.M. Cohen, D.B. Wales,
{\it Linearity of Artin groups of finite type},
preprint.

\smallskip
\item{[Cri]}
J. Crisp,
{\it Injective maps between Artin groups},
Cossey, John (ed.) et al., Geometric group theory down under.
Proceedings of a special year in geometric
group theory, Camberra, Australia, July 14--19, 1996.
De Gruyter, Berlin, 1999, pp. 119--137.

\smallskip
\item{[CrP]}
J. Crisp, L. Paris,
{\it The solution to a conjecture of Tits on the subgroup
generated by the squares of the generators of an Artin group},
Invent. Math., to appear.

\smallskip
\item{[Del]}
P. Deligne,
{\it Les immeubles des groupes de tresses g\'en\'eralis\'es},
Invent. Math. {\bf 17} (1972), 273--302.

\smallskip
\item{[Deo]}
V.V. Deodhar,
{\it On the root system of a Coxeter group},
Commun. Algebra {\bf 10} (1982), 611--630.

\smallskip
\item{[Dig]}
F. Digne,
{\it On the linearity of Artin braid groups},
preprint.

\smallskip
\item{[Hil]}
H. Hiller,
{\it Geometry of Coxeter groups},
Research Notes in Mathematics, 54. Pitman Advanced
Publishing program, Boston-London-Melbourn, 1982.

\smallskip
\item{[Kra1]}
D. Krammer,
{\it The braid group $B_4$ is linear},
Invent. Math. {\bf 142} (2000), 451--486.

\smallskip
\item{[Kra2]}
D. Krammer,
{\it Braid groups are linear},
preprint.

\smallskip
\item{[Mic]}
J. Michel, {\it A note on words in braid monoids},
J. Algebra {\bf 215} (1999), 366--377.

\smallskip
\item{[Tit1]}
J. Tits,
{\it Le probl\`eme des mots dans les groupes de Coxeter},
Sympos. Math., Roma 1, Teoria Gruppi, Dic. 1967 e Teoria 
Continui Polari, Aprile 1968, Academic Press, London,
1969, pp.175--185.

\smallskip
\item{[Tit2]}
J. Tits,
{\it Normalisateurs de tores. I: Groupes de Coxeter \'etendus},
J. Algebra {\bf 4} (1966), 96--116.

\bigskip\bigskip\noindent
\halign{#\hfill\cr
Luis Paris\cr
Laboratoire de Topologie\cr
Universit\'e de Bourgogne\cr
UMR 5584 du CNRS, BP 47870\cr
21078 Dijon cedex\cr
FRANCE\cr
\noalign{\smallskip}
{\tentt lparis@u-bourgogne.fr}\cr}

\end